\documentclass[11pt]{amsart}
\usepackage{graphicx}
\allowdisplaybreaks
\usepackage{amsmath,amssymb,mathrsfs,amsthm}

\usepackage{verbatim}
\usepackage[spanish,USenglish]{babel} % espanol, ingles
\usepackage{color}
\usepackage{tikz}
\usepackage{graphicx}

\newtheorem{thm}{Theorem}[section]
\newtheorem{corollary}[thm]{Corollary}
\newtheorem{proposition}[thm]{Proposition}
\newtheorem{theorem}[thm]{Theorem}
\newtheorem{lemma}[thm]{Lemma}

\theoremstyle{definition}
\newtheorem{definition}[thm]{Definition}

\theoremstyle{remark}
\newtheorem{remark}[thm]{Remark}

\numberwithin{equation}{section}

\def\N{\mathbb{N}}
\def\R{\mathbb{R}}

\def\C{\mathbb{C}}

\def\Z{\mathbb{Z}}

\def\Tn2{{\mathcal T}^{(n)}_2}
\def\LL{\mathcal{L}}

\def\NN{\mathbb{N}}
\def\RR{\mathbb{R}}

\def\CCma{\mathbb{C}^{+}}

\def\FF{\mathcal{F}}

\def\calR{\mathcal{R}}

\def\CC{\mathcal{C}}
\def\DD{\mathcal{D}}

\def\LL{\mathcal{L}}
\def\TT{\mathcal{T}}

\def\RR{\mathcal{R}}
\def\SS{\mathcal{S}}
\def\WW{\mathcal{W}}

\def\a{\alpha}

\def\ka{k_\alpha}

\def\G{\Gamma}

\def\LiiRma{L^{2}(\mathbb{R}^{+})}

\def\HpCma{H_{p}(\mathbb{C}^{+})}
\def\HqCma{H_{q}(\mathbb{C}^{+})}
\def\HrCma{H_{r}(\mathbb{C}^{+})}

\begin{document}

\title[RKH of Brownian type defined by Ces\`aro operators]
{RKH spaces of Brownian type defined by Ces\`aro-Hardy operators}

\author[Gal\'{e}]{Jos\'{e} E. Gal\'{e}}
\address{Departamento de Matem\'aticas, Instituto Universitario de Matem\'aticas y Aplicaciones, Universidad de Zaragoza, 50009 Zaragoza, Spain.}
\email{gale@unizar.es}

\author[Miana]{Pedro J. Miana}
\address{Departamento de Matem\'aticas, Instituto Universitario de Matem\'aticas y Aplicaciones, Universidad de Zaragoza, 50009 Zaragoza, Spain.}
\email{pjmiana@unizar.es}

\author[S\'{a}nchez-Lajusticia]{Luis S\'{a}nchez--Lajusticia}
\address{Departamento de Matem\'aticas, Facultad de Ciencias, Universidad de Zaragoza, 50009 Zaragoza, Spain.}
\email{luiss@unizar.es}

\thanks{Jos\'{e} E. Gal\'{e} and Pedro J. Miana have been partially supported by Project MTM2016-77710-P, and ID2019-105979GBI00, DGI-FEDER, of the MCEI and Project E48-20R, Gobierno de Arag\'on, Spain.}



\maketitle

\begin{abstract} We study reproducing kernel Hilbert spaces introduced as ranges of generalized Ces\`aro-Hardy operators, in one real variable and in one complex variable. Such spaces can be seen as formed by absolutely continuous functions on the positive half-line (or paths of infinite length) of fractional order, in the real case. A theorem of Paley-Wiener type is given which connects the real setting with the complex one. These spaces are related with fractional operations in the context of integrated Brownian processes. We give  estimates of the norms of the corresponding reproducing kernels.
\end{abstract}

%%%%%%%%%%%%%
\section{Introduction}
%%%%%%%%%%%%%

Let $\alpha>0$ and let
$f\colon\R^+\to\C$ be a measurable function. G. H. Hardy formulated \cite{H1} and proved in \cite{H2} the integral inequality
\begin{equation} \label{inicial}
\left(\int_0^\infty \left|\frac{1}{t^\a} \int_0^t (t-s)^{\a-1} f(s)ds\right|^p dt\right)^\frac{1}{p}
\leq M_\a\left(\int_0^\infty\vert f(t)\vert^p dt\right)^{1/p},
\end{equation}
for $1<p<\infty$ and optimal constant
$M_\a:=B(\a,1-{1\over p})=\frac{\G(\a)\G(1-\frac{1}{p})}{\G(\a +1-\frac{1}{p})}$
(see also \cite[Theorem 329]{HLP}), and
its dual inequality
\begin{equation}\label{dual}
\left(\int_{0}^{\infty} \left\vert\int_{t}^{\infty} \frac{(s-t)^{\a-1}}{s^\a} f(s)ds\right\vert^{p}dt\right)^\frac{1}{p}
\leq
M'_\a\left(\int_0^\infty\vert f(t)\vert^p dt\right)^{1/p},
\end{equation}
for $1\le p<\infty$ and optimal constant $M'_\a:=B(\a,{1\over p})$.  (Note that the inequality (\ref{inicial}) 
 also holds in the limit case $p=\infty$.) These estimates are generalizations of the case $\a=1$, which Hardy himself established as continuous versions of discrete inequalities that he had got in 1915, when he was searching for an elementary proof of Hilbert's inequality.
Early developments of Hardy's inequalities can be found in the classical book \cite{HLP}, and details about its history in both discrete and continuous forms are in \cite{KMP}, for instance.

The above inequalities induce linear bounded operators from $L_p(\R^+)$ into $L_p(\R^+)$ that we denote for $f\in L_p(\R^+)$ by
$$
\CC_\a(f)(t):={\a\over t^\a}\int_0^t (t-s)^{\a-1}f(s) ds, \hbox{ when } 1<p\le\infty,
$$
and
$$
\CC_\a^*(f)(t):=\a\int_t^\infty{(s-t)^{\a-1}\over s^\a}f(s) ds, \hbox{ when } 1\le p<\infty.
$$

For $\a=1$ operators $\CC_1$ or $\CC_1^*$, or their discrete counterparts, have received different names. As a sample, they are  called Hardy's operators in \cite{KMP}, \cite{DS}, Ces\`aro operators in    \cite{AS}, \cite{B}, \cite{BHS}, \cite{Mo1}, \cite {Mo2}, Copson operators in \cite{Mo1}, \cite {Mo2}, among other papers. There are also versions of the above operators in the complex plane  for $\alpha=1$, see \cite{AS}. 
The study of such operators is usually focused on problems around boundedness on diverse spaces, spectrum, interpolation, optimal domain and range,
\dots  (see for instance \cite{ABR}, \cite{AP}, \cite{DS}, \cite{LMPS}).
Here we call $\CC_\a$, $\CC_\a^*$ Ces\`aro-Hardy operators. We are interested in the range spaces, of these integral operators, endowed with the norm transferred from $L_p$ spaces, and more precisely in the Hilbertian case $L_2$. The motivation for such an approach is two-fold, arising from the connections of those operators with fractional integro-differentation, from one side, and  with fractional Brownian motion or white noise on the other hand.

Let $X$ be a Banach space and let $A\colon D(A)\to X$ a closed operator with domain $D(A)$ in $X$. In the study of the \lq\lq ill-posed" abstract Cauchy equation
$$
u'(t)=Au(t),\quad t\ge0;\, u(0)=x,\quad  x\in D(A)\subset X,
$$
so when the solution
$u\colon[0,\infty)\to X$
of the equation is not governed by a $C_0$-semigroup, families like $C$-semigroups or integrated semigroups, and homomorphisms like distribution semigroups are relevant, see \cite{ABHN}. In \cite{AK}, tempered distribution semigroups are considered which have as domains convolution Banach algebras -which we denote by
$\TT_1^{(n)}(t^n)$ here- defined, for $n\in\N$, as the completion of the space of test functions
$C_c^\infty(\R^+)$ in the norm
\begin{equation}\label{norm1}
\Vert f\Vert_{1,(n)}:=\int_0^\infty\vert f^{(n)}(t)\vert t^n\ dt<\infty, \quad f\in C_c^\infty(\R^+).
\end{equation}
(Similar algebras on the whole real line $\R$ had been introduced in \cite{BE}). The Banach algebra $\TT_1^{(n)}(t^n)$
admits an extension to fractional order of derivation $\a>0$ simply by considering certain fractional derivation (denoted by $W^\a f$) instead of the usual derivation $f^{(n)}$; see \cite{Mi} and \cite{GM}.
This extension, denoted by
$\TT_1^{(\a)}(t^\a)$, is also a convolution Banach algebra which has a number of applications related to functional calculi, integrated semigroups and theory of regular quasimultipliers, see \cite{GM}. Specific properties or applications of
$\TT_1^{(\a)}(t^\a)$ as a Banach algebra have been given in quite a number of papers, among them \cite{GMR1, GMR2, GMS, GS, GW}.
By replacing the $L_1$-norm of $t^{n}f^{(n)}$ with the $L_{p}$-norm, for $1<p \le \infty$, of $t^{n}f^{(n)}$ in (\ref{norm1}), one defines the convolution Banach
$\TT_1^{(\a)}(t^\a)$-module $\TT_p^{(\a)}(t^\a)$. It sounds sensible to find out properties and applications of such spaces similarly to the algebra case. A first analysis in that direction is done in \cite{GMMS} for $\a=n\in\N$ and $p=2$.

It turns out that the space $\TT_p^{(\a)}(t^\a)$ can be alternatively obtained as range space of the operator $\CC_\a^*$ with domain in 
$L_p(\R^+)$, so that  $\CC_\a^*$ may well be regarded under the viewpoint that fractional integro-differentiation provides.
On the other hand, fractional integrals and derivatives are of application in the theory of fractal Brownian motion
(fBm, for short) and self-similar systems (c.f., \cite{FP}, \cite{Hu}, \cite{M}, \cite{SL}), so that the
Ces\`aro-Hardy operators and the Hilbertian spaces that they define, namely $\TT_2^{(\a)}(t^\a)$, $\a>0$, appear in this way inserted in that theory. Furthermore, the action of the Laplace transform on $\TT_2^{(\a)}(t^\a)$ gives rise to a Hilbert space
$H_2^{(\a)}(\C^+)$ of holomorphic functions on the half-plane $\C^+:=\{z\in\C:\Re z>0\}$ which admits a simple description and might be of interest in order to deal with fBm of Riemann-Liouville type.

The organization of this paper is as follows.

For $t\in\R$, $1\le p\le\infty$ and $f\in L_p(\R^+)$, the formula $T_p(t)f(s):=e^{-t/p}f(e^{-t}s)$, for a.e. $s>0$, defines a $C_0$-group of isometries on $L_p(\R^+)$. 
 Subordination to $(T_p(t))_{t\in\R}$ or variations of  this (semi-)group has been considered in the case $\alpha=1$ to represent 
Ces\`aro-Hardy operators.For example, in relation with the subnormality of the Ces\`{a}ro operator \cite{C}, to study the Black-Scholes equation of financial mathematics \cite{AP} (the authors are grateful to V. Keyantuo for this observation), or to calculate norm and spectrum \cite{AP}, \cite{AS}. Subordination for $\alpha>0$, on the real line, is used in \cite{LMPS}.

We adopt this point of view  and use subordination to introduce Ces\`aro-Hardy operators
${\frak C}_\a$, ${\frak C}_\a^*$ and fractional integration-derivation of order $\alpha>0$ in the complex variable context
 (the expressions of ${\frak C}_1$, ${\frak C}_1^*$ under subordination to $T_p(t)$ are given in \cite[Th. 3.1 and Th. 3.3]{AS}). As a consequence of this procedure we easily show that the Laplace transform is an intertwining operator between real Ces\`aro-Hardy operators and complex Ces\`aro-Hardy operators. This is done in Section \ref{GCHO}.

In Section \ref{RangeCesar}, we deal with range spaces $\TT_p^{(\a)}(t^\a):=\CC_\a^*(L_p(\R^+))$,
for $\a>0$, endowed with the norm
$\Vert f\Vert_{p,(\a)}:=\Gamma(\a+1)\Vert (\CC_\a^*)^{-1}f\Vert_p$, for $f\in\TT_p^{(\a)}(t^\a)$. By expressing the operator $\CC_\a^*$ in terms of the Weyl operator $W^\a$ of the fractional calculus (see \cite{SKM} for the operators involved in the fractional calculus), the elements of  $\TT_p^{(\a)}(t^\a)$ can be seen as
absolutely continuous functions of fractional order $\a>0$. Proposition \ref{Tau-properties} contains several properties of the
space $\TT_p^{(\a)}(t^\a)$ and  its elements under that perspective. One of such properties is that point evaluations on
$\TT_p^{(\a)}(t^\a)$ have a sense and in fact are continuous whenever $\a>1/2$, whence in particular it follows that the space $\TT_2^{(\a)}(t^\a)$ is a reproducing kernel Hilbert space (RKHS, for short) for $\a>1/2$.

From Section \ref{RKH-Sobolev spaces} to the end of the paper we restrict our attention to the case $p=2$.
Section \ref{RKH-Sobolev spaces} is devoted to obtain the reproducing kernel which generates the Hilbert space $\TT_2^{(\a)}(t^\a)$, with $\a>1/2$, see Proposition \ref{kernelTaii}.
Then some ingredients of the general theory of RKHS are revised for $\TT_2^{(\a)}(t^\a)$, $\a>1/2$, from which it becomes apparent that there exists a relationship between this space and spaces arising in fractal Brownian motion or probability theory. This relation is partly described in Section \ref{brownian}, in connection with the Riemann-Liouville fractional calculus.

For $1\le p < \infty$, spaces $H_p^{(\a)}(\C^+)$ of holomorphic functions in $\C^+$, complex versions of spaces
$\TT_p^{(\a)}(t^\a)$,
are introduced in Section \ref{hardysobol}. To do this, one needs a complex form of the fractional calculus,
which is available through the subordination expression of the operator ${\frak C}_\a^*$ to the group $T_p(t)$.
Formally, just replacing real fractional derivatives with complex fractional derivatives, spaces $\TT_p^{(\a)}(t^\a)$
and $H_p^{(\a)}(\C^+)$ look identical. Moreover, for $p=2$ there is a correspondence of Paley-Wiener type in
the sense that
$\LL(\TT_2^{(\a)}(t^\a))=H_2^{(\a)}(\C^+)$ where $\LL$ is the Laplace transform. Indeed, 
 $H_2^{(\a)}(\C^+)$ is a
RKHS not only for
$\a>1/2$ but for all $\a>0$, and its reproducing kernel $K_\a$ can be expressed by means of a nice integral.
The Paley-Wiener result and the formula of the kernel are given in Theorem \ref{paleywiener}. This theorem
is a fairly wide extension of \cite[Th. 3.3 and Th. 4.1]{GMMS}, in the integer order case, to fractional order.
In Section \ref{estkernel}, it is proved that the function
$K_{\a,z}:=K_\a(\cdot, z)$ satisfies the estimate
$\Vert K_{\a,z}\Vert_{2,(\a)}\sim \vert z\vert^{-1/2}$, $z\in\C^+$, up to constants from below and from above.
Such an estimate is somehow surprising since, in
a number of classical examples
of kernels $\kappa\colon \Omega\times\Omega\to\C$ for spaces of holomorphic functions, usual estimates of the norm of $\kappa_y:=\kappa(\cdot,y)$
involve the distance of the point $y\in\Omega$ to the boundary of the domain $\Omega$, whereas
$\Vert K_{\a,z}\Vert_{2,(\a)}$
depends on the {\it radial} distance of $z$, that is, of
$z$ in $\C^+$ to the origin.

We have taken the operator $\CC_\a^*$ restricted on $L_2(\R^+)$ and on its range $\TT_2^{(\a)}(t^\a)$, as the way to
show the links of Ces\`aro-Hardy operators  with  fractional calculus.
This choice has been motivated by the fruitful relation of the spaces $\TT_2^{(\a)}(t^\a)$ with abstract Cauchy
equations and their associated families of bounded operators. Alternatively, we could have chosen to take the operator
$\CC_\a$ and range $\CC_\a(L_2(\R^+))$ and try to make a similar treatment.
The paper finishes with Section \ref{brownaverage}, which is a long remark rather than a section. There, we show that
$\TT_2^{(\a)}(t^\a)=\CC_\a(L_2(\R^+))$, which in view of the simple and good properties of spaces $\TT_2^{(\a)}(t^\a)$ and
$H_2^{(\a)}(\C^+)$ pointed out in the corresponding sections of this paper, suggests the question if averaging fractal operations, as $\CC_\a$ does, could be helpful within Brownian theory.

\medskip

%%%%%%%%%%%%%%%%%%
%%%%%%%%%%%%%%%%%%%
\section{Generalized Ces\`aro-Hardy operators}\label{GCHO}
%%%%%%%%%%%%%%%%%%%
%%%%%%%%%%%%%%%%%%%%

(i) {\it A group of isometries acting on the real line.} Let $f\colon\R^+\to\C$ be a complex, measurable function defined a.e. on the half-line $\R^+:=(0,\infty)$. For $t\in\R$ and $1\le p\le+\infty$, put
$$
T_p(t)f(s):=e^{-t/p}f(e^{-t}s),\quad \hbox { a.e. }\ s>0,
$$
where $t/\infty$ is understood with value $0$.
Clearly, $t\mapsto T_p(t)$ is a group in $t\in\R$ acting by composition on functions $f$ as above 
 with inverses $(T_p(t))^{-1}=T_p(-t)$, $t\in \R$. As a matter of fact,
$(T_p(t))_{t\in\R}$ is a strongly continuous group of  (surjective) isometries on $L_p(\R^+)$ if $1\le p<\infty$, and on
$C_0([0,\infty))$ if $p=\infty$ ($f\in C_0([0,\infty))$ means that $f$ is continuous in $[0,\infty)$ and $f(\infty)=0$). The isometric property of this group is fairly simple to check. As for the strong continuity, it is also part of folklore: for $1\le p<\infty$,
$s,t\in\R^+$ and
$h\in C_c(0,\infty)$ (this means that $h$ is continuous with compact support in $(0,\infty)$),
$$
\Vert T_p(t)h-T_p(s)h\Vert_{L_p(\R^+)}^p=\int_0^\infty\vert e^{-t/p}h(e^{-t}r)-e^{-s/p}h(e^{-s}r)\vert^p\ dr
$$
and so $\Vert T_p(t)h-T_p(s)h\Vert_{L_p(\R^+)}^p\to0$ as $t\to s$. For arbitrary $f$ in $L_p(\R^+)$, one obtains
$\Vert T_p(t)f-T_p(s)f\Vert_{L_p(\R^+)}^p\to0$, as $t\to s$, using the density of $C_c(0,\infty)$ and the fact that $T_p(t)$, $t\in\R$, are isometries. The case $p=\infty$ is even simpler.

(ii) {\it The same group (of isometries) acting on the complex plane.} Set $\C^+:=\{z\in\C:\Re z>0\}$. Let $\HpCma$ the Hardy space  on $\C^+$, which is formed by all holomorphic functions $F$ in $\C^+$ such that
$$
\Vert F\Vert_p:=\sup_{x>0}\left({1\over 2\pi}\int_{-\infty}^\infty\vert F(x+iy)\vert^p\ dy\right)^{1/p}<\infty.
$$
Endowed with the norm $\Vert \cdot\Vert_p$, $\HpCma$ is a Banach space. Recall that the elements $F$ of $\HpCma$ admit extension to the real line $i\R$ almost everywhere (via nontangential limits) with
$\Vert F\mid_{i\R}\Vert_{L_p(i\R)}=\Vert F\Vert_p$.

The space $\HpCma$ can alternatively be described as follows. For $\theta\in(-\pi/2,\pi/2)$ and every function $F$ holomorphic in $\C^+$, set $F_\theta(z):=F(ze^{i\theta})$, for $z\in\C^+$. Then $\HpCma$ is the space of holomorphic functions $F$ on $\C^+$ such that
$$
\Vert F\Vert_{p,rad}:=\sup_{-\pi/2<\theta<\pi/2}
\left({1\over2\pi}\int_0^\infty\vert F_\theta(r)\vert^pdr
\right)^{1/p}
$$
Moreover, $\Vert \cdot\Vert_{p,rad}$ is a norm and
\begin{equation}\label{radial}
\Vert F\Vert_{p,rad}=\Vert F\Vert_p, \quad \forall \ F\in \HpCma;
\end{equation}
see \cite{Se} (which extends to arbitrary $p$ the case $p=2$ proven in \cite{Dz}).

The operator $T_p(t)$ of part (i) above extends obviously to functions $F$ in the complex plane by putting
$T_p(t)F(z):=e^{-t/p}F(e^{-t}z)$, for all $z\in\C$, $t\in\R$. In particular, for $F\in\HpCma$, the function $T_p(t)F$ is holomorphic in $z\in\C^+$ and, for $1\le p<\infty$,
\begin{eqnarray*}
\Vert T_p(t)F\Vert_p^p&=&\sup_{x>0}\left({1\over 2\pi}\int_{-\infty}^\infty e^{-t}\vert F(e^{-t}(x+iy))\vert^p\ dy\right)\cr
&=&\sup_{x>0}\left({1\over 2\pi}\int_{-\infty}^\infty\vert F(e^{-t}x+iu)\vert^p\ du\right)=\Vert F\Vert_p^p,
\end{eqnarray*}
whence it follows that $T_p(t)$ is an isometry from $\HpCma$ onto itself for every $t\in\R$. Further,
\begin{eqnarray*}
\Vert T_p(t)F-T_p(s)F\Vert_p&=&\Vert(T_p(t)F)\mid_{i\R}-(T_p(s)F)\mid_{i\R}\Vert_p\cr
&=&\Vert T_p(t)(F\mid_{i\R})-T_p(s)(F\mid_{i\R})\Vert_p\to0, \hbox{ as } t\to s,
\end{eqnarray*}
similarly as above for the real case in part (i) (note that $T_p(t)$ is also an isometry on $L_p(i\R)$).
In conclusion, $(T_p(t))_{t\in\R}$ is a strongly continuous group of isometries on $\HpCma$, if $1\le p<\infty$.

When $p=\infty$ one needs to consider the Banach space ${\mathcal A}_0(\C^+)$ of all holomorphic functions $F$ in $\C^+$ with continuous extension to $\overline{\C^+}=\C^+\cup i\R$ and such that
$\lim_{z\in\overline{\C^+},\ z\to\infty}F(z)=0$. Then it is readily seen in a similar way to above that
$(T_\infty(t))_{t\in\R}$ is a strongly continuous group of isometries on ${\mathcal A}_0(\C^+)$.

(iii) {\it The group $T_p(t)$ and the Laplace transform.} Let $\mathcal L$ denote the usual Laplace transform for functions $f$ on
$\R^+$ such that the map $fe_z$, $z\in\C^+$, is integrable, where $e_z(t):=e^{-zt}$, $t>0$.
We have, just formally for the moment,
\begin{eqnarray*}
\LL(T_p(t)f)(z)&=&e^{-t/p}\int_0^\infty f(e^{-t}s)e^{-zs}\ ds=e^{t/q}\int_0^\infty f(r)e^{-ze^tr}\ dr\cr
&=&
e^{t/q}(\LL f)(e^tz)=T_q(-t)(\LL f)(z), \hbox{ for } t\in\R, z\in\C^+,
\end{eqnarray*}
where $1/q=(p-1)/p$ for $1\le p\leq\infty$. In short, $\LL\circ T_p(t)=T_q(-t)\circ\LL$.

In order to make the above relatonship rigorous, note that for every function $h\in C^{(2)}_c(\R^+)$ one has $h=e_1g$ with
$g\in C^{(2)}_c(\R^+)$ so that for every $z\in\C^+$,
$$
\vert \LL h(z)\vert=\vert(\LL g)(z+1)\vert={\vert(\LL g'')(z+1)\vert\over\vert(z+1)\vert^{2}}
\leq{(\LL(\vert g''\vert)(1)\over\vert(z+1)\vert^{2}}.
$$
Hence
\begin{eqnarray*}
\Vert(\LL h)\Vert_r
&\leq&
(\LL\vert g''\vert)(1)\sup_{x>0}\left({1\over2\pi}\int_{-\infty}^\infty{dy\over((x+1)^2+y^2)^r}\right)^{1/r}\cr
&=&(\LL\vert g''\vert)(1)\left({1\over2\pi}\int_{-\infty}^\infty{dy\over(1+y^2)^r}\right)^{1/r}<\infty;
\end{eqnarray*}
that is, {$\LL h \in\HrCma\cap {\mathcal A}_0(\C^+)$ for all  $1\le r  <\infty$}. Thus one has
$$
\LL(T_p(t)f)=T_q(-t)(\LL f), \quad t\in\R,\ \forall f\in C^{(2)}_c(\R^+).
$$

Certainly, the above equality holds with some restrictions for a wider class of functions. For example, let $p$ be such that $1\le p\le2$. Then Hausdorff-Young's inequality tells us that the Fourier transform
$\FF$ is bounded (it is a contraction, indeed) from $L_p(\R)$ into $L_q(\R)$. Now, for $x>0$ and $f\in L_p(\R^+)$ we have
$\LL f(x+i\ \cdot\ )=\FF(e_xf)\in L_q(\R)$ with
$\Vert\FF(e_xf)\Vert_{L_q(\R)}\le\Vert e_xf\Vert_{L_p(\R^+)}\le\Vert f\Vert_{L_p(\R^+)}$
(with the corresponding version when $q=\infty$). This means that $\LL f \in\HqCma$. So we have

\begin{equation}\label{laplacegrupo}
\LL(T_p(t)f)=T_q(-t)(\LL f), \quad \forall \ 1\le p\le2, f\in L_p(\R^+), t\in\R.
\end{equation}

%%%%%%%%%%%%%%%%%%%%%%%%%
\medskip
\subsection{Ces\`aro-Hardy operators on the real line}
%%%%%%%%%%%%%%%%%%%%%%%%%

Let $\alpha>0$. For a measurable function $f$ defined a.e. on $\R^+$ set
\begin{equation}\label{CesaroNoStar}
\CC_\alpha f(s):={\alpha\over s^\alpha}\int_0^s(s-u)^{\alpha-1}f(u)\ du, \quad s>0,
\end{equation}

\begin{equation}\label{CesaroStar}
\CC^*_\alpha f(s):=\alpha\int_s^\infty{(u-s)^{\alpha-1}\over s^\alpha}f(u)\ du, \quad s>0.
\end{equation}

Operators $\CC_\alpha$, $\CC^*_\alpha$ induced by formulas (\ref{CesaroNoStar}) and (\ref{CesaroStar}) will be called here (generalized) Ces\`aro-Hardy operators.
More precisely, take $f$ in $\SS_+:=\SS\mid_{[0,\infty)}$ where $\SS$ is the Schwarz class on
$\R$. Hardy's inequalities (\ref{inicial}), (\ref{dual}) say that
$$
\Vert\CC_\alpha f\Vert_{L_r(\R^+)}\le A_r\Vert f\Vert_{L_r(\R^+)}, \quad 1<r\le\infty,
$$
and
$$
\Vert\CC^*_\alpha f\Vert_{L_r(\R^+)}\le B_r\Vert f\Vert_{L_r(\R^+)}, \quad  1\le r<\infty,
$$
with constants
$A_r={\Gamma(\alpha+1)\Gamma(1-(1/r))\over\Gamma(\alpha+1-(1/r))}$ and
$B_r={\Gamma(\alpha+1)\Gamma(1/r)\over\Gamma(\alpha+(1/r))}$;
see \cite[p. 245]{HLP}.
Thus the integrals given in (\ref{CesaroNoStar})
and (\ref{CesaroStar}) define linear bounded operators  $\CC_\alpha\colon L_r(\R^+)\to  L_r(\R^+)$, $1<r\le\infty$, and
$\CC^*_\alpha\colon L_r(\R^+)\to  L_r(\R^+)$, $1\le r<\infty$, respectively.

In this paper, we initially focus on the operator $\CC^*_\alpha$ rather than on $\CC_\alpha$. Let $q,p$ be such that $1<q,p<\infty$ and $(1/q)+(1/p)=1$. It is readily seen, using Fubini's theorem, that
$\CC^*_\alpha\colon L_p(\R^+)\to  L_p(\R^+)$ is the
adjoint operator of $\CC_\alpha\colon L_q(\R^+)\to  L_q(\R^+)$. Moreover, the bounded operator
$\CC^*_\alpha\colon L_1(\R^+)\to  L_1(\R^+)$ can be regarded as the restriction operator to $L_1(\R^+)$ of the adjoint
$(\CC_\alpha\mid_{C_0(\R^+)})^*\colon{\mathcal M}(\R^+)\to {\mathcal M}(\R^+)$. Here ${\mathcal M}(\R^+)$ is the Banach space of bounded regular Borel measures on $\R^+$ and $\CC_\alpha\mid_{C_0(\R^+)}$ is the restriction of the bounded operator
$\CC_\alpha\colon L_\infty(\R^+)\to  L_\infty(\R^+)$ to $C_0(\R^+)$.

\medskip
$\bullet$ {\it Ces\`aro-Hardy operators as subordinated to semigroups:}

For $\alpha>0$, $1\le p<\infty$, $f\in L_p(\R^+)$ and $s>0$ we have 
\begin{eqnarray*}
\CC^*_\alpha f(s)
&=&\alpha\int_0^\infty(1-e^{-t})^{\alpha-1}e^{-t/p}(T_p(-t)f)(s)\ dt.
\end{eqnarray*}
 Similarly, we obtain for $1< q\le\infty$,
$g\in L_q(\R^+)$ and $s>0$,

\begin{eqnarray*}
\CC_\alpha f(s)
&=&\alpha\int_0^\infty(1-e^{-t})^{\alpha-1}e^{-t/p}(T_q(t)g)(s)\ dt
\end{eqnarray*}
where $p=q/(q-1)$ (\cite[Theorem 3.3 and 3.7]{LMPS}).

 Put $\varphi_{\a,p}(t):=\a(1-e^{-t})^{\alpha-1}e^{-t/p}$, $t>0$, for every $\alpha>0$ and $1\le p<\infty$.  Clearly 
$\varphi_{\a,p}\in L_1(\R^+)$,  so we get the following result.

\begin{proposition}\label{hillephil}
For $\alpha>0$, $1\le p<\infty$, $q$ such that $(1/p)+(1/q)=1$, $f\in L_p(\R^+)$ and  $g\in L_q(\R^+)$ we have
\begin{equation}\label{subordination}
\CC^*_\alpha f=\int_0^\infty\varphi_{\a,p}(t)T_p(-t)f\ dt,
\quad
\CC_\alpha g=\int_0^\infty\varphi_{\a,p}(t)T_q(t)g\ dt,
\end{equation}
where the convergence of integrals are in the Bochner sense.
\end{proposition}

The proposition tells us that $\CC^*_\alpha$ and $\CC_\alpha$ are given by subordination with respect to the semigroups
$(T_p(-t))_{t\ge0}$, $(T_q(t))_{t\ge0}$, respectively, in terms of the
Hille-Phillips operational calculus. This property has been observed for $\a=1$ in \cite{AS} and for $\a>0$ in \cite{LMPS} under slighty different expressions.

\begin{corollary}\label{conmuta}
For $\alpha, \beta>0$ and $1<p<\infty$,
$$
\CC_\a^*\circ\CC_\beta=\CC_\beta\circ\CC_\a^* \quad \hbox{ on }\, L_p(\R^+).
$$
\end{corollary}

\begin{proof}
Let $\alpha, \beta>0$, $1<p<\infty$ and $f\in L_p(\R^+)$. Since $(T_p(t))_{t\in\R}$ is a group,
\begin{eqnarray*}
(\CC_\a^*\circ\CC_\beta)f&=&\int_0^\infty\varphi_{\a,p}(t)T_p(-t)\left(\int_0^\infty\varphi_{\beta,q}(s)T_p(s)f\ ds\right)dt\cr
&=&\int_0^\infty\int_0^\infty\varphi_{\a,p}(t)\varphi_{\beta,q}(s)T_p(-t)T_p(s)f\ dsdt\cr
&=&\int_0^\infty\int_0^\infty\varphi_{\beta,q}(s)\varphi_{\a,p}(t)T_p(s)T_p(-t)f\ dsdt\cr
&=&\int_0^\infty\varphi_{\beta,q}(s)T_p(s)\left(\int_0^\infty\varphi_{\a,p}(t)T_p(-t)f\ dt\right)ds\cr
&=&(\CC_\beta\circ\CC_\a^*)f,
\end{eqnarray*}
as we wanted to show.
\end{proof}

Taking $\beta=\alpha$ in the corollary, one gets that $\CC_\a$ is normal on $L_2(\R^+)$. We do not deal with this property in this paper.

%%%%%%%%%%%%%%%%%%%%%%%%%%%%%%%%
\medskip
\subsection{Ces\`aro-Hardy operators on the (right-hand) complex half-plane}
%%%%%%%%%%%%%%%%%%%%%%%%%%%%%%%%%%%

We introduce here generalized Ces\`aro operators acting on $\C^+$ via the semigroup $T_p(t)$ of isometries considered in the preceding subsection. In this way, we avoid tedious calculations to check the holomorphy of the integral functions involved.

Let $1\le p<\infty$. Recall that $H_p(\C^+)$ is the Banach space of holomorphic functions $F$ on $\C^+$ such that
$$
\Vert F\Vert_p:=\sup_{x>0}\left({1\over 2\pi}\int_{-\infty}^\infty\vert F(x+iy)\vert^p\ dy\right)^{1/p}<\infty.
$$
In fact, $\Vert \cdot \Vert_p$ is the norm of $H_p(\C^+)$.

\begin{definition}\label{cesarcomplex}
For $\alpha>0$, $1\le p<\infty$, $q$ such that $(1/p)+(1/q)=1$, $F\in H_p(\C^+)$ and $G\in H_q(\C^+)$, define
$$
{\frak C}^*_\alpha F:
=\int_0^\infty\varphi_{\a,p}(t)T_p(-t)F dt\,
\hbox{ and }\,
{\frak C}_\alpha G:
=\int_0^\infty\varphi_{\a,p}(t)T_q(t)G dt,
$$
where the convergence of integrals are in the Bochner sense.
\end{definition}

From the definition it is evident that ${\frak C}^*_\alpha F$ is holomorphic in $\C^+$, in fact
${\frak C}^*_\alpha F\in H_p(\C^+)$ and so in particular there exists the nontangencial limit
${\frak C}^*_\alpha F(iy):=\lim_{z\to iy}{\frak C}^*_\alpha F(z)$ for almost everywhere $y\in\R$.

Let us develop the vector valued integral in Definition \ref{cesarcomplex}. Take $z\in\C^+$, $z=\vert z\vert e^{i\theta}$ with
$-\pi/2\le\theta\le\pi/2$ and $F\in H_p(\C^+)$. Put $F_\theta(u):=F(ue^{i\theta})$, for $u>0$. By (\ref{radial}) it follows that
$F_\theta\in L_p(\R^+)$ for all $\theta\in[-\pi/2,\pi/2]$. Then
\begin{eqnarray*}
{\frak C}^*_\alpha F(z)&=&\int_0^\infty\varphi_{\a,p}(t)T_p(-t)F(z) dt=\alpha\int_0^\infty(1-e^{-t})^{\alpha-1}F(ze^t) dt\cr
&=&\alpha\int_1^\infty u^{-\alpha}(u-1)^{\alpha-1}F_\theta(\vert z\vert u) du\cr
&=&\alpha\int_{\vert z\vert}^\infty(r-\vert z\vert)^{\alpha-1}r^{-\alpha}F_\theta(r) dr=\CC^*_\alpha F_\theta(\vert z\vert),
\end{eqnarray*}
where the latter operator $\CC^*_\alpha$ is the real one defined on $L_p(\R^+)$ as in (\ref{CesaroStar}).
Note that when
$\theta=-\pi/2$ or $\theta=\pi/2$ then ${\frak C}^*_\alpha F$ is defined only a. e. on $(-\infty,0)$ or $(0,\infty)$, respectively.

One can also express ${\frak C}^*_\alpha F$ on $\C^+$ by a complex integral:
For $z=\vert z\vert e^{i\theta}$ as above,
\begin{eqnarray*}
{\frak C}^*_\alpha F(z)&=&\alpha\int_{\vert z\vert}^\infty
(re^{i\theta}-\vert z\vert e^{i\theta})^{\alpha-1}(r e^{i\theta})^{-\alpha}F(re^{i\theta})d(re^{i\theta})\cr
&=&
\alpha\int_{\vert z\vert\cdot e^{i\theta}}^{\infty\cdot e^{i\theta}}
{(\lambda-z)^{\alpha-1}\over(\lambda)^\alpha}F(\lambda) d\lambda
=
\alpha\int_{0\cdot e^{i\theta}}^{\infty\cdot e^{i\theta}}
\lambda^{\alpha-1}{F(\lambda+z)\over(\lambda+z)^\alpha} d\lambda.
\end{eqnarray*}
Here $\lambda\mapsto \lambda^\beta$ is defined taking the principal argument continuous in
$\C\setminus (-\infty,0]$.

\begin{remark}
It can be shown by standard methods that the above complex integral giving  ${\frak C}^*_\alpha F$
on $\C^+$ is independent of the ray of integration; that is,
$$
{\frak C}^*_\alpha F(z)=\alpha\int_{0\cdot e^{i\omega}}^{\infty\cdot e^{i\omega}}
\lambda^{\alpha-1}{F(\lambda+z)\over(\lambda+z)^\alpha} d\lambda,
$$
for every $z\in\C^+$ and every $\omega\in(-\pi/2,\pi/2)$. We will not use this property here.

It can be also shown, with similar arguments as those used before, that one has
$$
{\frak C}_\a G(z)={\a\over z^\alpha}\int_{0\cdot e^{i\omega}}^{\vert z\vert e^{i\omega}}
(z-\lambda)^{\a-1}G(\lambda) d\lambda,
$$
for $G$, $z$ and $\omega$ as above, after Definition \ref{cesarcomplex}. We do notice this fact about the operator
${\frak C}_\a$ for the sake of completeness, but it will not be used in this paper.
\end{remark}

An immediate consequence of Definition \ref{cesarcomplex} is that the complex Ces\`aro-Hardy operators commute.

\begin{corollary}\label{complexcesarocommu}
For $\a,\beta>0$, $1<p<\infty$ and $F\in H_p(\C^+)$,
$$
({\frak C}_\a^*\circ{\frak C}_\beta)F=({\frak C}_\beta\circ{\frak C}_\a^*)F.
$$
\end{corollary}

Another interesting consequence of the subordination to the groups $T_r(t)$ is that
the Laplace transform $\LL$ intertwines Ces\`aro-Hardy operators.

\begin{corollary}\label{laplacesaro}
Let $\a>0$. Then,
\item{\rm (1)} \quad $\LL\circ\CC^*_\a=\frak C_\a\circ\LL \, \hbox{ on } L_p(\R^+) \hbox{ for } 1\le p\le2$.
\item{\rm (2)} \quad $\LL\circ\CC_\a=\frak C_\a^*\circ\LL \, \hbox{ on } L_q(\R^+) \hbox{ for } 1< q\le2$.
\end{corollary}

\begin{proof}
For $1\le p\le2$, $z\in\C^+$ and $f\in L_p(\R^+)$ one has
\begin{eqnarray*}
[\LL(\CC^*_\alpha f)](z)&=&
\int_0^\infty\varphi_{\a,p}(t)[(\LL\circ T_p(-t))f(z)\ dt\cr
&=&
\int_0^\infty\varphi_{\a,p}(t)[(T_q(t)\circ\LL)f](z)\ dt=[\frak C_\alpha(\LL f)](z),
\end{eqnarray*}
where we have used (\ref{laplacegrupo}) in the second equality. This gives us part (1). Part (2) is shown analogously.
The proof is over.
\end{proof}

%%%%%%%%%%%%%%%%%%%%%%%%
%%%%%%%%%%%%%%%%%%%%%%%%
\medskip
\section{Range spaces of Ces\`aro-Hardy operators}\label{RangeCesar}
%%%%%%%%%%%%%%555
%%%%%%%%%%%%%%%%

The optimal domain and/or optimal range of distinguished operators have been studied in recent times. The case of Ces\`aro operators is discussed in \cite{DS}.   
 Here we are interested in the range of operators $\CC^*_\alpha$ when the domains are $L_p$ spaces as above. Actually, such range spaces can be identified with spaces $\TT_p^{(\alpha)}(t^\alpha)$ introduced using the Weyl fractional derivation operator $W^\alpha$, see \cite{GM, LMPS, Mi, MRS}:

 To begin with, recall that
$\CC^*_\alpha\colon L_p(\R^+)\to L_p(\R^+)$, $1\le p<\infty$, is injective. We give a proof of this for the convenience of readers.

Let $p$ be such that $1\le p<\infty$ and take $f\in L_p(\R^+)$ such that $\CC^*_\alpha f=0$. Then, for every $t>0$,
$$
\CC^*_\alpha f(t)=\alpha\int_t^\infty{(s-t)^{\alpha-1}\over s^\alpha}f(s)ds
=\alpha t^{\alpha-1}\int_0^{1/t}(t^{-1}-r)^{\alpha-1}r^{-1}f(r^{-1})dr.
$$
Therefore $\CC^*_\alpha f= 0 \Leftrightarrow((\cdot)_+^{\alpha-1}\ast g)=0$ where $g(r):=r^{-1}f(r^{-1})$. In addition,
$g\in L_{loc}^1(\R^+)$ because
$\int_0^A g(r)dr=\int_{1/A}^\infty{f(t)\over t}dt$ for every $A>0$. By Titchmarsh's convolution theorem \cite[p. 188]{D}, one has $g=0$ and so $f=0$.
(The injectivity of $\CC^*_\alpha$ is also a consequence of the density of the range of  $\CC_\alpha$ on $L_q$ spaces but we have preferred to recall the above direct argument.)

Then we define $\TT_p^{(\alpha)}(t^\alpha):= \CC^*_\alpha(L_p(\R^+))$ endowed with the norm
$$
\Vert f\Vert_{p,(\alpha)}:=\Gamma(\a+1)\Vert (\CC^*_\alpha)^{-1}f\Vert_{L_p(\R^+)},
$$
so that $\TT_p^{(\alpha)}(t^\alpha)$ is a Banach space and $\CC^*_\alpha\colon L_p(\R^+)\to \TT_p^{(\alpha)}(t^\alpha)$ is an (onto) isometry.

Spaces $\TT_p^{(\alpha)}(t^\alpha)$, $\alpha\ge0$ -and therefore the Ces\`aro operators $\CC^*_\alpha$, $\CC_\alpha$- are intimately related with fractional derivatives and integrals:

Let $L_p(\R^+,t^{\a p})$ denote the Banach space of measurable functions $f$ such that the function $t\mapsto t^{\alpha}f(t)$ belongs to $L_p(\R^+)$, and let
$\tau_\alpha$ denote the multiplication operator by the (weight) function $t\mapsto t^{\alpha}$, $t>0$. Put
$\mu_{-\a}:=\Gamma(\alpha+1)\tau_{-\a}$.

Set $W^\alpha\colon{\mathcal T}_p^{(\alpha)}(t^\alpha)
\buildrel { ({\mathcal C}_\alpha^*)^{-1}  }
\over\longrightarrow L_p(\R^+)\buildrel{\mu_{-\a}}\over\longrightarrow L_p(\R^+,t^{\alpha p})$;
that is,
$$
W^\alpha f(t):=\Gamma(\alpha+1)t^{-\alpha}\left[({\mathcal C}_\alpha^*)^{-1}f\right](t),
\quad f\in\mathcal T_p^{(\alpha)}(t^\alpha), t>0.
$$
It is clear that $W^\alpha$ has the inverse
$$
W^{-\alpha} g(t):=\int_t^\infty(s-t)^{\alpha-1} g(s) {ds\over\Gamma(\alpha)}, \quad g\in L_p(\R^+,t^{\alpha p}).
$$

Thus $\TT_p^{(\alpha)}(t^\alpha)$ is formed by all the elements $f$ in $L_p(\R^+)$ for which there exists a unique element in $L_p(\R^+,t^{\alpha p})$,
notated $W^\alpha f$, such that
\begin{equation}\label{repres}
f(t):=\int_t^\infty(s-t)^{\alpha-1} W^{\alpha} f(s) {ds\over\Gamma(\alpha)},
\end{equation}

On account of the above notation the expression of the norm of $f$ in $\TT_p^{(\alpha)}(t^\alpha)$ has the form
$$
\Vert f\Vert_{p,(\alpha)}:=
\left(\int_0^\infty\vert W^\alpha f(t)t^\alpha\vert^pdt\right)^{1/p}<\infty.
$$

Notice that in the limit, when $\a=0$, we get $\TT_p^{(0)}(t^0)=L_p(\R^+)$. The operators $W^{\alpha}$ and $W^{-\alpha}$ introduced above are extensions of the corresponding restricted operators
$W^{\alpha}\colon\SS_+\to\SS_+$, $W^{-\alpha}\colon\SS_+\to\SS_+$ which can be found
in \cite{SKM}, \cite{MR}, for instance, as operators defining Weyl fractional derivation
and Weyl integration, respectively. The following properties emphasize the derivation character of $W^\alpha$.
Set $W^0:=Id$, the identity operator, and $h_\lambda(t):=h(\lambda t)$, for any function
$h$ and $\lambda, t>0$.

\begin{proposition}\label{Wderiv}
\item{\rm(i)} {\rm Integro-differentiation group property:} $W^\alpha\circ W^\beta=W^{\alpha+\beta}$
on $\SS_+$ for every $\alpha,\beta\in\R$.
\item{\rm(ii)} $W^n\varphi=(-1)^n\varphi^{(n)}$, for every $\varphi\in \SS_+$, if $\alpha=n\in\Z$. Hence, for every $\alpha>0$
and every integer $n$ such that $n\ge[\alpha]+1$,
$$
W^\alpha\varphi=(-1)^n{d^n\over dt^n}W^{-(n-\alpha)}\varphi.
$$
\item{\rm(iii)} {\rm Homogeneity:} $W^\alpha f_\lambda=\lambda^\alpha(W^\alpha f)_\lambda$, where
$f\in\TT_p^{(\alpha)}(t^\alpha)$.
\end{proposition}

\begin{proof}
(i) See \cite[p. 96]{SKM}.

(ii) For a negative integer $n$, the first equality is the formula of $-n$ times integration by parts. For positive $n$, the equality follows from the equality $W^n=(W^{-n})^{-1}$.

(iii) This equality is straightforward.
\end{proof}

In fact, for later considerations, it is suitable to regard spaces $\TT_p^{(\alpha)}(t^\alpha)$ as being formed by
\lq\lq derivatives" of functions. Under this viewpoint, these spaces were introduced in \cite{Mi}, \cite{GM} (and previously in \cite{AK} in the case
$\alpha\in\N$)  for $p=1$ in the setting of integrated semigroups and distribution semigroups (these families of \lq\lq semigroups" are of interest to deal with ill-posed abstract Cauchy problems, see \cite{ABHN}). The case $p>1$ was introduced in \cite{LMPS, MRS, R}. Next, we list some of their properties.

\begin{proposition}\label{Tau-properties}
Let $1\le p<\infty$.
\item{\rm(i)} $C_c^{(\infty)}(\R^+)$ is dense in $\TT_p^{(\alpha)}(t^\alpha)$ for all $\alpha\ge0$.
\item{\rm(ii)} For every $\alpha>0$, $\TT_p^{(\alpha)}(t^\alpha)$ is a convolution Banach $\TT_1^{(\alpha)}(t^\alpha)$-module;
that is, there exists a constant $M_{\alpha,p}$ such that for every $g\in\TT_1^{(\alpha)}(t^\alpha)$ and
$f\in\TT_p^{(\alpha)}(t^\alpha)$,
$$
g\ast f\in\TT_p^{(\alpha)}(t^\alpha)
\hbox{ with } \Vert g\ast f\Vert_{p,(\a)}\le M_{\alpha,p}\Vert g\Vert_{1,(\a)} \Vert f\Vert_{p,(\a)}.
$$
Moreover, $\TT_1^{(\alpha)}(t^\alpha)\ast\TT_p^{(\alpha)}(t^\alpha)$ is dense in $\TT_p^{(\alpha)}(t^\alpha)$.
\item{\rm(iii)} For every $\alpha, \beta$ such that $\beta>\alpha\ge0$,
$$
{\TT}_p^{(\beta)}(t^\beta)\hookrightarrow {\TT}_p^{(\alpha)}(t^\alpha)
\hookrightarrow L_p(\R^+).
$$
(the hook arrows mean continuous inclusions).

Furthermore, for every $f\in\TT_p^{(\beta)}(t^\beta)$,
\begin{equation}\label{escala}
W^\alpha f(t)={1\over\Gamma(\beta-\alpha)}\int_t^\infty(s-t)^{\beta-\alpha-1}W^\beta f(s)ds, \quad t>0,
\end{equation}
whence
\begin{equation}\label{Winequal}
\vert W^\alpha f(t)\vert\le C_{\alpha,\beta,p} t^{-(\alpha+(1/p))}\Vert f\Vert_{p,(\alpha)}, \quad t>0, \hbox{ provided } \beta>\alpha+{1\over p},
\end{equation}
for some constant $C_{\alpha,\beta,p}>0$.
\item{\rm(iv)} For every $\alpha>0$,
\begin{eqnarray*}
{\TT}_p^{(\alpha)}(t^\alpha)&=&\{f\in L_p(\R^+): t^\alpha W^\alpha f\in L_p(\R^+)\}\cr
&=&\{f\in L_p(\R^+): t^\nu W^\nu f\in L_p(\R^+),\ \forall\ 0\le\nu\le\alpha\}
\end{eqnarray*}
and the following norms are equivalent on ${\TT}_p^{(\alpha)}(t^\alpha)$:
$$
(1)\ \Vert f\Vert_{p,(\alpha)}:=\Vert t^\alpha W^\alpha f\Vert_{L_p(\R^+)},\quad
$$
$$
(2)\ \sup_{0\le\nu\le\alpha}\Vert f\Vert_{p,(\nu)}
=\sup_{0\le\nu\le\alpha}\Vert t^\nu W^\nu f\Vert_{L_p(\R^+)}.
$$
\end{proposition}

\begin{proof}
To begin with, we do observe two facts. First, $C_c^{(\infty)}([0,\infty))$ is dense in every
$\TT_p^{(\alpha)}(t^{\alpha})$
for all $\alpha>0$ and $1\le p<\infty$, since
$\CC_\alpha^*\colon L^p(\R^+)\to\TT_p^{(\alpha)}(t^{\alpha})$ is a surjective isometry with
$\CC_\alpha^*(C_c^{(\infty)}([0,\infty)))=C_c^{(\infty)}([0,\infty))$.
 On the other hand, take $\varphi\in C_{c}^{(\infty)}(\R^+)$ positive such that $\int_{0}^{\infty}\varphi(t)\,dt =1$.
For $\varepsilon >0$ put
$\varphi_\varepsilon(t) = \varepsilon^{-1}\varphi(\varepsilon^{-1}t)$, $t \in \R^{+}$.
Then $(\varphi_\varepsilon)_{0<\varepsilon<1}$ is a bounded approximate identity for $\TT_p^{(\alpha)}(t^{\alpha})$
for every $\alpha\ge0$, that is,
$\lim_{\varepsilon\to0^+} f\ast\varphi_\varepsilon=f$
in $\TT_p^{(\alpha)}(t^{\alpha})$ (see \cite[Proposition 2.3]{GMR2} for $p=1$; for arbitrary $p$ the argument is similar).

(i) Since $h\ast\varphi_\varepsilon\in C_{c}^{(\infty)}(\R^+)$ for every $h\in C_{c}^{(\infty)}([0,\infty))$ and
$C_c^{(\infty)}([0,\infty))$ is dense in $\TT_p^{(\alpha)}(t^{\alpha})$, one gets that $C_c^{(\infty)}(\R^+)$ is dense in
$\TT_p^{(\alpha)}(t^{\alpha})$.

(ii) Since $C_c^{(\infty)}(\R^+)$ is dense in $\TT_p^{(\alpha)}(t^{\alpha})$ one can derive the module property of the statement as it is done in \cite{R}.

(iii) Let $\alpha>0$ and take $\beta$ such that $\beta>\alpha\ge0$. Let $f\in\TT_p^{(\beta)}(t^{\beta})$. The function
$\phi$ given by the integral
$$
\phi(t):=\int_t^\infty(s-t)^{\beta-\alpha-1}\vert W^\beta f(s)\vert ds, \quad t>0,
$$
is an element of $L_p(t^{\alpha p})$, so that that integral is finite for all a.e. $t>0$. In effect,
\begin{eqnarray*}
\Vert\phi\Vert_{L_p(t^{\a p})}&\le&\left(\int_0^\infty\left[\int_t^\infty{(s-t)^{\beta-\alpha-1}\over s^{\beta-\alpha}}
s^\beta \vert W^\beta f(s)\vert ds\right]^pdt\right)^{1/p}\cr
&\le&M_{\alpha, \beta}\Vert s^\beta W^\beta f(s)\Vert_{L_p(\R^+)}=M_{\alpha, \beta}\Vert f\Vert_{p,(\beta)}<\infty,
\end{eqnarray*}
for some positive constant $M_{\alpha, \beta}$, where the second inequality is Hardy's inequality (\ref{dual}).

Hence, the function $g$ given by
$$
g(t):={1\over\Gamma(\beta-\alpha)}\int_t^\infty(s-t)^{\beta-\alpha-1}W^\beta f(s)\ ds,
$$
is defined for a.e. $t>0$ and $g\in L_p(t^{\a p})$ with
$\Vert g\Vert_{L_p(t^{\a p})}\le M'_{\alpha, \beta}\Vert f\Vert_{p,(\beta)}$, for some constant $M'_{\alpha, \beta}$.

Using the same argument as before
but with the function $\psi$ given by $\psi(t):=\int_t^\infty(s-t)^{\alpha-1}\phi(s)\ ds<\infty$, for $t>0$, instead of $\phi$ (so in particular
with $\alpha-1$ in the exponent instead
$\beta-\alpha-1$) one obtains that $\psi\in L_p(\R^+)$ whence  $\psi(t)<\infty$ for every $t>0$.

\bigskip

Then by Fubini's theorem and (\ref{repres}), for every $t>0$,
\begin{eqnarray*}
\int_t^\infty {(s-t)^{\alpha-1} \over\Gamma(\alpha)} g(s)\ ds & = & \int_t^\infty\int_t^r {(s-t)^{\alpha-1}\over\Gamma(\alpha)} {(r-s)^{\beta-\alpha-1} \over \Gamma(\beta-\alpha)}  ds \ W^\beta f(r)\ dr \cr
& = & \int_t^\infty { (r-t)^{\beta-1} \over\Gamma(\beta)} W^\beta f(r)\ dr = f(t).
\end{eqnarray*}
Therefore, applying the uniqueness of the representation (\ref{repres}) we
derive that $g=W^\alpha f$.
In other words, we have proved (\ref{escala}) and the continuous inclusion
${\TT}_p^{(\beta)}(t^\beta)\hookrightarrow {\TT}_p^{(\alpha)}(t^\alpha)$.

Finally, take $\beta>\alpha+(1/p)$ and $f\in\TT_p^{(\beta)}(t^{\beta})$.

Assume $1< p<\infty$. For $t>0$,
\begin{eqnarray*}
\vert W^\alpha f(t)\vert&\le&{1\over\Gamma(\beta-\alpha)}
\int_t^\infty{(s-t)^{\beta-\alpha-1}\over s^\beta} s^\beta\vert W^\beta f(t)\vert\ ds\cr
&\le&
{1\over\Gamma(\beta-\alpha)}
\left(\int_t^\infty(s-t)^{q(\beta-\alpha-1)}s^{-q\beta}ds\right)^{1/q}\Vert f\Vert_{p,(\beta)}\cr
&=&
{B(q(\alpha+1)-1,q(\beta-\alpha-1)+1)^{1/q}\over\Gamma(\beta-\alpha)}
t^{-(\alpha+{1\over p})}\Vert f\Vert_{p,(\beta)}.
\end{eqnarray*}

Assume now $p=1$. Then, for $t>0$,
\begin{eqnarray*}
\vert W^\alpha f(t)\vert&\le&{1\over\Gamma(\beta-\alpha)}
\left(\sup_{s>t}{(s-t)^{\beta-\alpha-1}\over s^\beta}\right)\Vert f\Vert_{1,(\beta)}\cr
&=&M_{\alpha,\beta}t^{-(\alpha+1)}\Vert f\Vert_{1,(\beta)}.
\end{eqnarray*}

(iv) This point follows readily from the definiton of $\TT_p^{(\a)}(t^\a)$, equality (\ref{escala})
and Hardy's inequality (\ref{dual}).

All in all, the proof is over.
\end{proof}

\begin{remark}\label{AbsCont}
\normalfont
For $n\in\N$, it follows from Proposition \ref{Wderiv} and (\ref{Winequal}) that every function $f$ in ${\TT}_p^{(n)}(t^n)$ is ($n-1$)-times differentiable with $f^{(n-1)}$ absolutely continuous on $\R^+$ and such that
$\int_0^\infty\vert f^{(n)}(t)\vert^pdt<\infty$. This suggests us to refer to ${\TT}_p^{(\alpha)}(t^\alpha)$ as space of absolutely continuous functions of fractional order, when $\alpha$ is any positive number. In order to distinguish the class
$\TT_p^{(\alpha)}(t^\alpha)$, $\alpha>0, 1\le p<\infty$, from other classes of Sobolev type existing in the literature,
each $\TT_p^{(\alpha)}(t^\alpha)$ will be called Lebesgue-Sobolev space here.
\end{remark}

Part (ii) of Proposition \ref{Tau-properties} says in particular that ${\TT}_1^{(\alpha)}(t^\alpha)$ is a Banach algebra.
In this respect, ${\TT}_1^{(\alpha)}(t^\alpha)$ has been studied in a number of papers (see \cite{GMR1}, \cite{GMR2},  \cite{GS}, and references therein, for example).
It sounds sensible to also study the structure of spaces ${\TT}_p^{(\alpha)}(t^\alpha)$, $1< p<\infty$.
From (\ref{Winequal}), it follows that $\vert f(t)\vert\le C\Vert f\Vert_{p,(\alpha)}t^{-1/p}$, for every $\alpha>1/p$,
$f\in{\TT}_p^{(\alpha)}(t^\alpha)$ and $t>0$. This is to say that point evaluations are continuous on
${\TT}_p^{(\alpha)}(t^\alpha)$. When $p=2$ it means that spaces ${\TT}_2^{(\alpha)}(t^\alpha)$  are Hilbert spaces with reproducing kernel. We call them RKH-Sobolev spaces and show, in the sequel of the present paper, some of their features.

%%%%%%%%%%%%%%%%%%%%%%%
%%%%%%%%%%%%%%%%%%%%%%%
\section{RKH-Sobolev spaces}\label{RKH-Sobolev spaces}
%%%%%%%%%%%%%%%%%%%%%
%%%%%%%%%%%%%%%%%%%%%

In this section and in the remainder of the paper, we focus on the case $p=2$. It is clear that the space
${\TT}_2^{(\alpha)}(t^\alpha)$ is, for every $\alpha>0$, a Hilbert space with inner product
$$
(f|g)_{2,(\alpha)}:=\int_{0}^{\infty}W^{\alpha}f(t)\overline{W^{\alpha}g(t)}t^{2\alpha}dt,
\quad f,g \in{\TT}_2^{(\alpha)}(t^\alpha).
$$

Via the isometry constructed out from $(\CC^*_\alpha)^{-1}$, one can find a suitable ortonormal basis in
${\TT}_2^{(\alpha)}(t^\alpha)$.  Let $(\ell_{m})_{m=0}^{\infty}$ be the orthonormal system on $L_2(\R^+)$ of Laguerre
functions $\ell_m$ given by
$$
\ell_{m}(t) = e^{-t/2} \sum_{j=0}^{m} \binom{m}{j} \frac{(-1)^{j}}{j!}t^{j}, \quad t>0, m=0,1,\dots
$$
Set $\ell_{m,\alpha}:=W^{-\alpha}(t^{-\alpha}\ell_{m})$, that is, for $t>0$,
\begin{equation}\label{OrtoBasis}
\ell_{m,\alpha}(t)
=\sum_{j=0}^{m} \binom{m}{j} \frac{(-1)^{j}}{j!} \frac{1}{\Gamma(\alpha)}
\int_{1}^{\infty} (u-1)^{\alpha-1} u^{j-\alpha}u^{j}e^{-ut/2} du.
\end{equation}
Then $(\ell_{m,\alpha})_{m=0}^{\infty}$ is an orthonormal basis in ${\TT}_2^{(\alpha)}(t^\alpha)$ since $W^{-\alpha}$ is an isometry from $L_2(t^{2\alpha})$ onto ${\TT}_2^{(\alpha)}(t^\alpha)$. This basis will be used in Section \ref{hardysobol}.

As it has been pointed out at the end of Section \ref{RangeCesar}, for $f\in{\TT}_2^{(\alpha)}(t^\alpha)$ one has that
$f(t)$ exists and
$\vert f(t)\vert\le t^{-1/2}C\Vert f\Vert_{2,(\alpha)}$ for every $t>0$ and $\alpha>1/2$. Thus,
${\TT}_2^{(\alpha)}(t^\alpha)$ is a reproducing kernel Hilbert space (RKHS for short). Our aim next is to find an expression of the reproducing kernel in ${\TT}_2^{(\alpha)}(t^\alpha)$.

Let $t>0$. The space ${\TT}_2^{(\alpha)}(t^\alpha)$ is a RKHS if and only if
there exists $k_{\alpha,t}\in{\TT}_2^{(\alpha)}(t^\alpha)$ such that $(f\mid k_{\alpha,t})_{2,(\alpha)}=f(t)$
for all
$f\in{\TT}_2^{(\alpha)}(t^\alpha)$. On the other hand, $f(t)=
{1\over\Gamma(\alpha)}\int_t^\infty (u-t)^{\alpha-1}W^\alpha f(u)du
={1\over\Gamma(\alpha)}\int_0^\infty W^\alpha f(u) (u-t)_+^{\alpha-1}du$, and therefore
$$
\int_0^\infty W^\alpha f(u)\overline {W^\alpha k_{\alpha,t}(u)}u^{2\alpha}du=
\int_0^\infty W^\alpha f(u) (u-t)_+^{\alpha-1}{du\over\Gamma(\alpha)}
$$
for every $f$ in ${\TT}_2^{(\alpha)}(t^\alpha)$, so for every $f\in C_c^{(\infty)}(\R^+)$. It follows then that
\begin{equation}\label{Wkernel}
W^\alpha k_{\alpha,t}(u)={1\over\Gamma(\alpha)}{(u-t)_+^{\alpha-1}\over u^{2\alpha}}, \quad t,u>0.
\end{equation}
Note that the last function is in $L_2(u^{2\a})$ if and only if $\alpha>1/2$. In this case we have
\begin{eqnarray*}
k_{\alpha,t}(s) &=& W^{-\alpha}\left[ W^{\alpha}k_{\alpha,t}\right](s)
= \frac{1}{\Gamma(\alpha)} \int_{s}^{\infty} (r-s)^{\alpha - 1} W^{\alpha}k_{\alpha,t}(u) du\cr
&=&
\frac{1}{\Gamma(\alpha)^2} \int_{0}^{\infty} (r-s)_{+}^{\alpha- 1}(r-t)^{\alpha - 1}_+ r^{-2\alpha} dr, \quad s>0.
\end{eqnarray*}
In conclusion, we have proved the following result. Put $k_\alpha(s,t):=k_{\alpha,t}(s)$.

\begin{proposition}\label{kernelTaii}
Let $\alpha>0$. The Hilbert space ${\TT}_2^{(\alpha)}(t^\alpha)$ is RKHS if and only if $\alpha>1/2$. In this case, the kernel for
${\TT}_2^{(\alpha)}(t^\alpha)$ is the function
\begin{equation}\label{kernelT}
k_\alpha(s,t) = \int_{0}^{\infty} g_\alpha(s,r) g_\alpha(t,r) dr, \quad t,s>0,
\end{equation}
where
$$
g_\alpha(t,r) = \frac{(r-t)^{\alpha - 1}_{+}}{r^\alpha \Gamma(\alpha)}, \quad t,r>0.
$$
\end{proposition}

\begin{remark}
\normalfont
Spaces $\TT_2^{(\a)}(t^\a)$ can be considered as spaces formed by paths (of infinite lenght in this case), as it happens with other typical examples in the theory of reproducing kernels.
One can readily describe some standard or general facts of such a theory in our setting.

(1) {\bf Norm of the kernel}. For $\alpha>1/2$,
\begin{eqnarray*}
\Vert k_{\alpha,t}\Vert_{2,(\alpha)}^2&=&k_\alpha(t,t)={1\over\Gamma(\alpha)^2}
\int_t^\infty(s-t)^{2\alpha-2}s^{-2\alpha}ds\cr
&=&{1\over\Gamma(\alpha)^2}
\int_t^\infty\left(1-{t\over s}\right)^{2\alpha-2}s^{-2}ds={1\over\Gamma(\alpha)^2t}
\int_0^1(1-u)^{2\alpha-2}du\cr
&=&
{B(1,2\alpha-1)\over\Gamma(\alpha)^2}{1\over t}
={1\over\Gamma(\alpha)^2 (2\alpha-1)}{1\over t},
\end{eqnarray*}
where $B$ is the beta function. So,
$$
\Vert k_{\alpha,t}\Vert_{2,(\alpha)}={1\over\Gamma(\alpha)\sqrt{2\alpha-1}}{1\over\sqrt {t}}, \quad t>0.
$$

(2) {\bf Hypergeometric function}. The kernel $k_{\alpha}$ can be  rewriten in terms of the hypergeometric function $_{2}F_{1}$. Recall that
for $\Re(c) > \Re(b) > 0$ and $|z|<1$, that function can be expressed by
$$
B(b,c-b) _{2}F_{1}(a,b,c,z) = \int_{0}^{1}t^{b-1}(1-t)^{c-b-1}(1-zt)^{-a}dt.
$$
Denote $t\wedge s =\min(t,s)$ and $t \vee s=\max(t,s)$.
Then
\begin{eqnarray*}
k_{\alpha}(s,t) & = & \int_{0}^{\infty} \frac{(r-s)^{\alpha - 1}_{+}}{r^{\alpha} \Gamma(\alpha)}
\frac{(r-t)^{\alpha - 1}_{+}}{r^{\alpha} \Gamma(\alpha)} dr \\
& = & \int_{t \vee s}^{\infty} \frac{r^{\alpha-1} \left(1-\frac{s}{r}\right)^{\alpha-1} r^{\alpha-1}
\left(1-\frac{t}{r}\right)^{\alpha-1}}{\Gamma(\alpha)^{2} r^{2\alpha}} dr \\
& = & \int_{0}^{\frac{1}{s}\wedge\frac{1}{t}} \frac{(1-su)^{\alpha-1}(1-tu)^{\alpha-1}}{\Gamma(\alpha)^{2}} du \\
& = & \int_{0}^{1}\frac{\left(1-\left(1\wedge\frac{s}{t}\right)y\right)^{\alpha-1}\left(1-\left(1\wedge\frac{t}{s}\right)y\right)^{\alpha-1}}{(s\vee t)\Gamma(\nu)^{2}} dy \\
& = & \frac{1}{(s\vee t)\Gamma(\alpha)\Gamma(\alpha+1)} {}_{2}F_{1}\left(1-\alpha, 1, \alpha+1, \frac{s \wedge t}{s \vee t} \right). \\
\end{eqnarray*}
When $\alpha=n\in\N$ one obtains
\begin{eqnarray*}
k_n(s,t) & = & \frac{1}{(s\vee t)(n-1)! n!} \sum_{j=0}^{n-1} (-1)^{j}\binom{n-1}{j}\frac{(1)_{j}}{(n+1)_{j}}\left(\frac{s\wedge t}{s\vee t}\right)^{j}\cr
& = &
 \sum_{j=0}^{n-1} \frac{(-1)^{j}}{(n+j)!(n-j-1)!} \frac{(s\wedge t)^{j}}{(s\vee t)^{j+1}},
\end{eqnarray*}
where $(a)_{b}$ is the Pochhammer symbol, and the second equality is a simplification.

If $\alpha=1$, we get   $k_1(s,t)=\frac{1}{s \vee t} = \frac{1}{s} \wedge \frac{1}{t}={\frak b}_{0}({1\over s},{1\over t}) $, where ${\frak b}_{0}(s,t) := s \wedge t$ is the well known reproducing kernel of the RKHS related to Brownian motion (or the covariance of the Brownian process). Moreover, the preceding  expressions of the kernel $k_\a$ resembles kernels arising in fractional versions of the Brownian motion.
This takes us to find out the relationship between spaces $\TT^{(\alpha)}_{2}(t^{\alpha})$ and spaces of the Brownian motion.
We will consider this item in Section \ref{brownian}.
\end{remark}

(3) {\bf Green function.}
Given a differential operator and  the corresponding equation $Lu(x) = f(x)$, the Green function for $L$ is the solution
(whenever it exists) $v=G$
to the twin equation
$Lv(x,s)=\delta(x-s)$, where $\delta$ is the Dirac delta distribution. Then, once $G$ has been found, one obtains
$u(x) = \int G(x,s)f(s)ds$ as solution to the initial equation. From the reproducing kernel theory and from (\ref{Wkernel}),
we have that $g_{\alpha}$ is Green's function for a certain operator $L_g$. This operator is given in terms of the fractional differential operator
$W^\alpha$ or the Ces\`aro operator. In effect,
\begin{eqnarray*}
g_\alpha(t,r) &=& \frac{(r-t)^{\alpha - 1}_{+}}{r^\alpha \Gamma(\alpha)}
\Longleftrightarrow
r^\alpha g_\alpha(t,r) = \frac{(r-t)^{\alpha - 1}_{+}}{ \Gamma(\alpha)}\cr
&\Longleftrightarrow&W^{\alpha}\left(s^{\alpha} g_{\alpha}(\cdot,s)\right)(t) = \delta_{s}(t),
\end{eqnarray*}
so that
$$
L _h= W^{\alpha} \circ \mu_{-\alpha}=  \mu_{\alpha}\circ (\CC^*_\alpha)^{-1} \circ \mu_{-\alpha}.
$$

(4) {\bf Green kernel integral.}
It is well known that the Green function asociated with a reproducing kernel allows us to recover the Hilbert space that the kernel generates via an integral transformation. In our case such a transformation ${\frak T}$ is, up to constants, the Ces\`aro-Hardy operator $\CC_\a^*$:
$$
{\frak T} f(t) = \int_{0}^{\infty}g_{\alpha}(t,r)f(r)dr
=\int_t^\infty\frac{(r-t)^{\alpha - 1}}{\Gamma(\alpha)r^\alpha}f(r)dr,\quad t>0, f\in L_2(\R^+);
$$
that is, ${\frak T}=\Gamma(\alpha+1)^{-1}\CC_\alpha^*$ and so $\TT^{(\alpha)}_{2}(t^\alpha)={\frak T}(L_2(\R^+))$,
as it had to be, see  \cite[Th. 11.3, Cor. 11.4]{PR} or \cite[Th. 1, p.4]{S2}.

%%%%%%%%%%%%%%%%%%%%%%%
%%%%%%%%%%%%%%%%%%%%%%%
\section{RKH-Sobolev spaces and Brownian motion}\label{brownian}
%%%%%%%%%%%%%%%%%%%%%
%%%%%%%%%%%%%%%%%%%%%

Given a definite-positive kernel $k$, there exists a (unique) Gaussian, zero mean, stochastic process $X_{t}$ such that the covariance is given by the kernel, which is to say $Cov(X_{t},X_{s})=k(t,s)$.
Let $B_{t}$ denote the well known Brownian motion, or Wiener process, whose covariance is
$$
{\frak b}_0(s,t)=\min\{s,t\}= \int_{0}^{t\wedge s} \chi_{(0,s)}(u)\chi_{(0,t)}(u) du, \quad s,t>0.
$$
 There are many textbooks where to find the Brownian motion and its main properties; see \cite{Du} for example. With the aim to provide useful models for the study of random phenomena with a strong interdependence between distant samples, the Brownian motion was widen to fractional Brownian motion (fBm, for short) by B. B. Mandelbrot and J. W. Van Ness in their seminal paper \cite{MV}. The $n$-times integrated Brownian motion $B_{n,t}$  can be defined recursively, in terms of stochastic integration, see \cite{L} and references therein. For $n=0$ we retrieve the Brownian motion, and the case $n=1$ is usually called the Langevin process. The order of integration need not be a positive  integer; 
in fact, P. L\`{e}vy had already introduced in \cite{Le} the Holmgren-Riemann-Liouville fractional integral $B_{\a,t}$ of order 
$\a>0$ of $B_t$ (with covariance ${\frak b}_\alpha$, say). 

Another important stochastic process is the so-called white noise, which can be considered as the formal derivative of the Brownian motion (despite the Brownian motion is nowhere differentiable with probability one, see \cite[p.140]{P}).
The corresponding integrated white noise is $N_{\a,t}=B_{\a-1,t}$, see \cite{SL}.

In practice, the kernel (covariance) of the fBm is fairly complicated, and in fact fBm is not suitable for modelling phenomena only arising in positive time, for example. Thus fBm is modified to simplify computations or to deal with specific type of problems, see \cite{BA}, \cite{FP}, \cite{Hu}, \cite{SL} for instance. In the above references, as well as in many other works, the tool used to approach questions involving Brownian phenomena is the fractional (integral) calculus. Several function spaces  have been considered as appropriate models to work out problems in fractal Brownian motion, see \cite{FP}, \cite{Hu} as samples of pertinent references. Note that also the space $\TT_2^{(\a)}(t^\a):=\CC^*_\a(L_2(\R^+))
=W^{-\a}(L_2(\R^+,t^{2\a}))$ lies in   this setting \cite[p. 5]{Hu}. Next, we introduce a model space $\RR^{(\alpha)}_2$ connected with Riemann-Liouville fractional calculus.

For $\a>0$, let $t^{\a-1}_{+}$ denote the function $t^{\a-1}_{+}:=t^{\a-1}\chi_{(0,\infty)(t)}$, $t\in\R$, and
let $\DD^{-\a}$ be the Riemann-Liouville integral operator
$$
\DD^{-\a}f(x) := \int_{0}^{x} \frac{(x-y)^{\a-1}}{\G(\a)} f(y) dy, \qquad x>0, \quad f\in\LiiRma.
$$

Note that the above integral is finite for a.e. $x>0$ since
$\displaystyle\DD^{-\a}f(x)=x^\a\left({1\over x^\a}\int_{0}^{x} \frac{(x-y)^{\a-1}}{\G(\a)} f(y) dy\right)$, so
$\DD^{-\a}f\in\tau^\a L_2(\R^+)$ by Hardy's inequality (\ref{inicial}).
Since $\DD^{-\a}f=t^{\a-1}_{+}\ast f$, Titchmarsh's convolution theorem implies that $\DD^{-\a}$ is injective.
Define $\RR^{(\alpha)}_2:=\DD^{-\a}(L_2(\R^+))$, so that for all $\varphi\in\RR^{(\alpha)}_2$ there exists a unique
$\DD^{\a}\varphi$ in $L_2(\R^+)$ 
such that
$$
\varphi(x)=\int_{0}^{x} \frac{(x-y)^{\a-1}}{\G(\a)} \DD^{\a}\varphi(y) dy, \quad x>0.
$$
 The function $\DD^{\a}\varphi$ is the well-known Riemann-Liouville fractional derivative of $\varphi$. We endow 
$\RR^{(\alpha)}_2$ with the norm $\|\varphi\|_{\RR_{(\alpha),2}}:=\|\DD^{\a}\varphi\|_{L^{2}}$.

\begin{remark} The space $\RR^{(1)}_2$ can be found for example in
\cite[p. 149]{PR} in the form $\RR^{(1)}_2= \{f : [0,+\infty)\rightarrow\CC \ | \ f \textrm{ absolutely continuous, } f(0)=0,  f'\in L^{2} \}$, or in \cite[p. 14]{SS} as $\RR^{(1)}_2= \{f \in W^{1,2}(0,\infty) \ | \ \lim_{\varepsilon\downarrow 0}f(\varepsilon)=0 \}$, where $W^{1,2}(0,\infty)$ is the Sobolev space of differential order 1 based on $L^{2}$. According to \cite[ p. 243]{BT},
$\RR^{(1)}_2$ is called the Cameron-Martin space and its unit ball known as the Strassen set. For $\a=n\in\NN$,
one has that
$\RR^{(n)}_2$ is isometrically isomorphic to the subspace of the Sobolev space $W^{n,2}(0,1)$ formed by
the functions $f$ satisfying the boundary conditions $f^{(j)}(0)=0$, $j=0,\cdots,n-1$,  see \cite[p. 92]{BT}.
Variants of the space $\RR^{(\alpha)}_2$, like $\DD^{-\a}(L_2([0,1])$ for instance, have been considered in the literaure, see  \cite{FP}, \cite{Hu}.
\end{remark}

\medskip
The covariance or kernel associated with the space
$\RR^{(\alpha)}_2$ is given  by the formula
$$
{\frak n}_{\a}(t,s) = \int_{0}^{t\wedge s} \frac{(t-u)^{\a-1}}{\G(\a)}  \frac{(s-u)^{\a-1}}{\G(\a)} du,
$$
with $\a>0$ referring to the number of \lq\lq times'' that we integrate the white noise (see \cite{FP}, \cite{SL}).
(${\frak b}_{\a}={\frak n}_{\a+1}$.)

We also have
$$
{\frak n}_{\a}(t,s)=\frac{(ts)^{\a -1}(t \wedge s)}{\G(\a)\G(\a+1)} {}_{2}F_{1}\left(1-\a, 1, \a+1, \frac{t \wedge s}{t \vee s} \right)
$$
whence
$$
{\frak n}_{\a}(t,s) = (t\vee s)(t\wedge s) (ts)^{\a-1} \ka(t,s) = (ts)^{\a} \ka(t,s).
$$

We next show a natural isometry between $\TT_2^{(\a)}(t^\a)$ and $\RR_2^{(\a)}$. This isometry is likely part of the folklore, but we have been unable to find a place where to get it explicitly.

\begin{lemma} \label{derivF}
Let $f\in\TT_2^{(\alpha)}(t^\alpha)$ and define $\varphi(x)=x^{\a-1}f\left(\frac{1}{x}\right)$ for $x>0$. Then
$$
\varphi\in\RR_2^{(\alpha)}\quad \hbox{ and }\quad \DD^{\a}\varphi(x) = x^{-(\a+1)}W^{\a}f\left(\frac{1}{x}\right), \quad x>0.
$$
\end{lemma}

\begin{proof}
For $f,\varphi$ as in the statement,
\begin{eqnarray*}
\varphi(x)= x^{\a-1}f\left(\frac{1}{x}\right) \Leftrightarrow \varphi(x)
&=& \frac{x^{\a-1}}{\G(\a)} \int_{1/x}^{\infty} \left(t-\frac{1}{x} \right)^{\a-1} W^{\a} f(t)dt\cr
&=& \frac{1}{\G(\a)} \int_{0}^{x} (x-y)^{\a-1} y^{-(\a+1)} W^{\a} f\left(\frac{1}{y}\right)dy,
\end{eqnarray*}
with
$$
\int_{0}^{\infty} \! \Big| y^{-(\a+1)} W^{\a} f\left(\frac{1}{y}\right) \Big|^{2} dy = \int_{0}^{\infty} \!\!|t^{\a} W^{\a} f(t)|^{2} dt < \infty.
$$

Thus the proof is over.
\end{proof}

Note that the mapping
$$
\begin{array}{cccc} L^{2}(0,\infty) & \longrightarrow & L^{2}(0,\infty) \\
                            f(x)      &   \mapsto   &  x^{-1}f(x^{-1})  \\
\end{array}
$$
is an isometric isomorphism.
Here is the isometry we referred to.

\begin{proposition} \label{isoTaHBa}
The mapping defined by
$$
\begin{array}{ccccc} \Theta_{\a} : & \TT_2^{(\alpha)}(t^\alpha) & \longrightarrow & \RR_2^{(\alpha)} \\
                                  &      f    &   \mapsto   &  x^{\a-1}f(x^{-1})\\
\end{array}
$$
is an isometric isomorphism.
\end{proposition}

\begin{proof} By Lemma \ref{derivF}, $\Theta_{\a}$ is well defined; moreover, it is obviously injective. As for the surjectivity, note that if
$\varphi\in\RR_2^{(\alpha)}$ then $\DD^\a\varphi\in L_2(\R^+)$,
that is, $\frac{1}{x}\DD^\a\varphi\left(\frac{1}{x}\right)\in L_2(\R^+)$ as noticed prior to the proposition. Therefore,
there exists a unique $f\in\TT_2^{(\alpha)}(t^\alpha)$ such that
$\frac{1}{x}\DD^\a\varphi\left(\frac{1}{x}\right) = x^{\a} W^{\a} f(x)$ and we know from Lemma \ref{derivF} that
$$
\DD^\a\varphi(x) = x^{-(\a+1)} W^{\a} f\left(\frac{1}{x}\right)  \Leftrightarrow \varphi(x) = x^{\a-1}f\left(\frac{1}{x}\right),
$$
as we wanted to show.
\end{proof}

As a consequence of the proposition one finds the following continuous inclusion.

\begin{corollary}
For every $\a>0$,
$$
\RR^{(\alpha)}_2\hookrightarrow L^{2}(\R^+,x^{-2\a}).
$$
\end{corollary}

\begin{proof}
For $\varphi\in\RR^{(\alpha)}_2$,
\begin{eqnarray*}
\int_{0}^{\infty} |\varphi(x)|^{2} \frac{dx}{x^{2\alpha}} &=&
\int_{0}^{\infty} \Big|x^{\a-1}f\left(\frac{1}{x}\right)\Big|^{2} \frac{dx}{x^{2\a}}\cr
&=& \int_{0}^{\infty} \Big| \frac{1}{x} f\left(\frac{1}{x}\right)\Big|^{2} dx = \displaystyle\int_{0}^{\infty} |f(x)|^{2} dx < \infty
\end{eqnarray*}
for a (unique) $f\in\TT_2^{(\alpha)}(t^\alpha)$. Moreover, for some constant $M_\a$,
\begin{eqnarray*}
\Vert \varphi\Vert_{L_2(x^{-2\a})} &=& \Vert f\Vert_2
\leq M_\a\Vert f\Vert_{2,(\alpha)}
=M_\a\int_{0}^{\infty} \vert t^{\a}W^{\a}f(t)\vert^{2} dt\cr
&=&M_\a\int_{0}^{\infty} \Big| x^{-(\a+1)} W^{\a}f\left(\frac{1}{x}\right)\Big|^{2} dx\cr
&=&M_\a \|\DD^\a\varphi\|_2
=M_\a \|\varphi\|^2_{\RR^{(\alpha)}_2}.
\end{eqnarray*}
\end{proof}

\begin{remark}
By Proposition \ref{Tau-properties} (iii), $\TT_2^{(\mu)}(t^\mu) \hookrightarrow \TT_2^{(\nu)}(t^\nu)$ for $\mu>\nu$. In contrast, there is no continuous embeddings between the spaces $\calR^{(\mu)}_2, \calR^{(\nu)}_2$: Let $F(t)= t/(1+t)$. Then $F\notin\LiiRma$, but $F'(t)=1/(1+t)^2 \in\LiiRma$, and $F\in\calR^{(1)}_2$.
\end{remark}

%The Laplace transform of functions in $\RR^{(\a)}_2$ is easy to obtain. 
For $r\in\R$, put $\zeta_r(z):=z^r=e^{r\log z}$,
$z\in\C^+$, where $\log z$ is the principal branch of the logarithm with principal argument in $[-\pi,\pi)$.  The following is the Paley-Wiener type result corresponding to
the Hilbert space $\RR^{(\a)}_2$.

\begin{corollary}\label{tLforR}
For all $\varphi\in\RR^{(\a)}_2$ and $z\in\C^+$,
$$
\LL(\DD^\a\varphi)(z) = z^{\a}\LL(\varphi)(z) \quad \hbox{and therefore} \quad \LL(\RR^{(\a)}_2)
=\zeta_{-\a} H_2(\C^+).
$$
\end{corollary}

\begin{proof}
For $\varphi\in\RR^{(\a)}_2$ we have $\varphi=\G(\a)^{-1}(t^{\a-1}_+\ast\DD^\a\varphi)$. Hence,
$$
\LL(\varphi)=\LL\left({t^{\a-1}_+\over\G(\a)}\ast\DD^\a\varphi\right)
=\LL\left({t^{\a-1}_+\over\G(\a)}\right)\LL(\DD^\a\varphi)=\zeta_{-\a}\LL(\DD^\a\varphi),
$$
which proves the first equality. Thus we have $\LL(\RR^{(\a)}_2)\subseteq\zeta_{-\a} H_2(\C^+)$.

Conversely, suppose that $F$ is a holomorphic function in $\C^+$ such that $\zeta_\a F\in H_2(\C^+)$. Then, by Paley-Wiener's theorem, there is $\phi\in L_2(\R^+)$ such that $\LL(\phi)=\zeta_\a F$. Hence,
$F=\zeta_{-\a}\LL(\phi)=\LL(\G(\a)^{-1}t^{\a-1}_+)\LL(\phi)=\LL(\G(\a)^{-1}t^{\a-1}_+\ast\phi)=\LL(\varphi)$ with
$\varphi=\G(\a)^{-1}t^{\a-1}_+\ast\phi\in\RR^{(\a)}_2$ and the proof is over.
\end{proof}

By Corollary \ref{tLforR}, it is readily seen that the space $\LL(\RR^{(\a)}_2)$, endowed with the norm
$\vert\Vert F\vert\Vert:=\Vert\zeta_\a F\Vert_2$, $F\in\zeta_{-\a} H_2(\C^+)$, is a RKHS with reproducing kernel given by
$(z,w)\in\C^+\times\C^+\mapsto(z\overline{w})^{-\a}(z+\overline w)^{-1}\in\C$. 
We deal with the space $\LL(\TT^{(\a)}_2(t^\a))$ and the representation of its reproducing kernel in the next section.

%%%%%%%%%%%%%%%%%%%%%%%
%%%%%%%%%%%%%%%%%%%%%%%
\section{Hardy-Sobolev spaces}\label{hardysobol}
%%%%%%%%%%%%%%%%%%%%%
%%%%%%%%%%%%%%%%%%%%%

Recall, for $1\le p<\infty$ and $F\in H_p(\C^+)$,
$$
T_p(t)F(z):=e^{-t/p}F(e^{-t}z), \quad t\in\R, z\in\C^+,
$$
is a $C_0$-group of isometries on $H_p(\C^+)$.

Let ${\frak C}_\alpha^*$ be the Ces\`aro-Hardy operator, introduced in Definition \ref{cesarcomplex}, given by
$$
{\frak C}^*_\alpha F:=\int_0^\infty\varphi_{\a,p}(t)T_p(-t)F\ dt\in H_p(\C^+),
$$
where $\varphi_{\a,p}(t):=\a(1-e^{-t})^{\alpha-1}e^{-t/p}$;  $\alpha>0$, $t>0$.
As seen after that definition,
${\frak C}_\alpha^*F(z)={\CC}_\alpha^*F_\theta(\vert z\vert)$ for $F\in H_2(\C^+)$, $z\in\C^+$, whence it follows that
${\frak C}_\alpha^*$ is injective.

\begin{definition}\label{HSespacio}
For $1\le p<\infty$, define the Hardy-Sobolev space, of order $\alpha>0$, $H_p^{(\alpha)}(\C^+)$ by
$$
H_p^{(\alpha)}(\C^+):={\frak C}_\alpha^*(H_p(\C^+)),
$$
endowed with the norm
$\Vert F\Vert_{p,(\alpha)}:=\Gamma(\a+1)\Vert ({\frak C}_\alpha^*)^{-1}F\Vert_p$, $F\in H_p^{(\alpha)}(\C^+)$.
\end{definition}

Then, in analogy to the real case, put
$$
\WW^\alpha:=\Gamma(\a+1)[\zeta_{-\alpha}\circ ({\frak C}_\alpha^*)^{-1}]
$$

or, equivalently,

$$
{\WW}^{-\alpha}={1\over\Gamma(\alpha+1)}{\frak C}_\alpha^*\circ \zeta_{\alpha},
$$
where $\zeta_{\nu}$ is the mutiplication operator by $z^\nu$, for $\nu\in\R$.

With this operational notation, we have that $F\in H^{(\alpha)}_p(\C^+)$ if and only if there exists
$\WW^\alpha F$ holomorphic in $\C^+$, with $\zeta_\alpha \WW^\alpha F\in H_p(\C^+)$, such that
\begin{equation}\label{weylcomplex}
F(z)={1\over\Gamma(\alpha)}\int_z^\infty(\lambda-z)^{\alpha-1}\WW^\alpha F(\lambda)d\lambda,\ z\in\C^+,
\end{equation}
where the integration path is the ray connecting $z$ with the complex infinity point. Also,
$\Vert F\Vert_{p,(\alpha)}:=\Vert \zeta_\alpha W^\alpha F\Vert_p$.

Indeed, working on rays leaving the origin in $\overline\C^+$, one gets, like in (\ref{escala}),
$$
\WW^\nu F(z)
={1\over\Gamma(\alpha-\nu)}\int_z^\infty(\lambda-z)^{\alpha-\nu-1}\WW^\alpha F(\lambda)d\lambda,\ z\in\C^+,
$$
for every $\nu$ such that $0\leq\nu<\alpha$. From here, one obtains the continuous inclusions
$$
H_p^{(\beta)}(\C^+)\hookrightarrow H_p^{(\alpha)}(\C^+)\hookrightarrow H_p^{(0)}(\C^+)=H_p(\C^+),
$$
for all $\beta>\alpha$.

From now on, we consider the case $p=2$. We know that
$H_2(\C^+)$ is a RKHS with kernel $K(z,w)=(z+\overline{w})^{-1}$; $z,w\in\C^+$.
Hence
one has that point evaluations
$$
\hbox{ev}_z\colon H^{(\alpha)}_2(\C^+)\hookrightarrow H_2(\C^+)\to\C,\quad z\in\C^+,
$$
are continuous on $H^{(\alpha)}_2(\C^+)$ and so this space is a RKHS. Let $K_\alpha$ denote its reproducing kernel,
so that $K_{\alpha,w}$ belongs to $H^{(\alpha)}_2(\C^+)$, where
$$
K_{\alpha,w}(z):=K_{\alpha}(z,w),\quad z,w\in\C^+.
$$

The aim of this section is to establish the following theorem. Its first part provides an integral expression for the kernel
$K_\alpha$. The second part is a Paley-Wiener type result.
\begin{theorem}\label{paleywiener} Let $\alpha>0$.
\item{\rm(i)} The space $H^{(\alpha)}_2(\C^+)$ is a RKHS with reproducing kernel
$$
K_\a(z,w)= \int_{0}^{1} \int_{0}^{1} \frac{(1-x)^{\a-1}}{\G(\a)} \frac{(1-y)^{\a-1}}{\G(\a)} \frac{1}{xz+y\overline w} dx dy,\ z,w\in\C^+.
$$
\item{\rm(ii)} The Laplace transform is an isometric isomorphism from $\TT_2^{(\alpha)}(t^\alpha)$ onto $H^{(\alpha)}_2(\C^+)$.
\end{theorem}

We prove the theorem using the basis method.
Let $(\ell_m)_{m\ge0}$ be the orthonormal basis of Laguerre functions in $L_2(\R^+)$  and let
$\ell_{m,\alpha}:=W^{-\alpha}(t^{-\alpha}\ell_{m})$ be the orthonormal basis in ${\TT}_2^{(\alpha)}(t^\alpha)$ obtained from
$(\ell_m)_{m\ge0}$ via the isometry
$W^{-\alpha}(t^{-\alpha}(\cdot))$, both bases given in the beginning of
Section \ref{RKH-Sobolev spaces}.

\begin{lemma}\label{laplacelaguerre}
For $\a>0$, $m\in\N\cup\{0\}$ and $z\in\C^+$,
$$
\leqno(i)\qquad\qquad (\LL\ell_{m})(z)=\displaystyle\frac{2(2z-1)^{m}}{(2z+1)^{m+1}}\, ,
$$
and
$$
\leqno(ii)\qquad \qquad \LL(\ell_{m,\a})(z) = \displaystyle
\frac{2}{\G(\a)}\int_{1}^{\infty}  \frac{(u-1)^{\a-1}}{u^{\a}} \frac{(2z-u)^{m}}{(2z+u)^{m+1}} du.
$$
\end{lemma}

\begin{proof}
For $m\in\N\cup\{0\}$ and $z\in\C^+$,
\begin{eqnarray*}
(\LL\ell_{m})(z) &=& \int_{0}^{\infty} e^{-zu} e^{-u/2}\sum_{j=0}^{m} \binom{m}{j} (-1)^{j} \frac{u^{j}}{j!} du\cr
&=& \sum_{j=0}^{m} (-1)^{j} \binom{m}{j} \frac{1}{j!} \frac{j!}{\left( z+\frac{1}{2}\right)^{j+1}}
= {2\over 2z+1}\sum_{j=0}^{m} \binom{m}{j} \left(-\frac{2}{2z+1} \right)^{j}\cr
&=&\frac{2}{2z+1} \left(1-\frac{2}{2z+1} \right)^{m} = \frac{2(2z-1)^{m}}{(2z+1)^{m+1}}.
\end{eqnarray*}

Also,
\begin{eqnarray*}
\LL(\ell_{m,\a})(z)& = & \sum_{j=0}^{m} \binom{m}{j} \frac{(-1)^{j}}{j!} \int_{1}^{\infty}
\int_{0}^{\infty} t^{j} e^{-t\left( \frac{u}{2}+z\right)}dt {(u-1)^{\a-1}\over u^{\a-j}} \frac{du}{\G(\a)}\cr
&=&\frac{1}{\G(\a)} \int_{1}^{\infty}  \frac{(u-1)^{\a-1}}{u^{\a}}
\sum_{j=0}^{m} \binom{m}{j} \left(\frac{-u}{\frac{u}{2}+z} \right)^{j}\frac{du}{\frac{u}{2}+z}\cr
&=& \frac{1}{\G(\a)} \int_{1}^{\infty}  \frac{2(u-1)^{\alpha-1}}{u^{\alpha}(u+2z)} \left(1-\frac{2u}{u+2z} \right)^{m}du\cr
& = & \frac{2}{\G(\a)}\int_{1}^{\infty} \frac{(u-1)^{\alpha-1}}{u^{\a}} \frac{(2z-u)^{m}}{(2z+u)^{m+1}} du,
\end{eqnarray*}
and we wanted to show.
\end{proof}

Since $(\LL(\ell_m))_{m\ge0}$ is an orthonormal basis in $H_2(\C^+)$ and the mapping
$$
H_2(\C^+)\to H_2^{(\alpha)}(\C^+),\, F\mapsto \WW^{-\a}(\zeta_{-\a}F)
$$
is an isometry, it follows that
$$
{\frak L}_{m,\a}:=\WW^{-\a}(\zeta_{-\a}\LL\ell_m), \, m=0,1,\dots ,
$$
form an orthonormal basis in $H_2^{(\alpha)}(\C^+)$.

\begin{lemma} \label{baseH2a}
For $\a>0$, $m\in\N\cup\{0\}$ and $z\in\C^+$,
\begin{eqnarray*}
{\frak L}_{m,\a}(z) &=& \frac{2}{\Gamma(\a)}\int_{1}^{\infty}\frac{(u-1)^{\a-1}}{u^{\a}}\frac{(2uz-1)^{m}}{(2uz+1)^{m+1}}du\cr
&=&{(-1)^m\over2z}\LL(\ell_{m,\a})({1\over4z}).
\end{eqnarray*}
\end{lemma}

\begin{proof} Let $z=\vert z\vert e^{i\theta}\in\CCma$, $m\in\NN\cup\{0\}$. Then, integrating on the ray (path) given by $\theta$
and applying Lemma \ref{laplacelaguerre} (i),
\begin{eqnarray*}
{\frak L}_{m,\a}(z) &=&\frac{2}{\Gamma(\alpha)}
\int_{z}^{\infty} \frac{(\lambda-z)^{\alpha-1}}{\lambda^{\alpha}} \frac{(2\lambda-1)^{m}}{(2\lambda+1)^{m+1}} d\lambda\cr
&=&\frac{2}{\Gamma(\alpha)} \int_{z}^{\infty} z^{\alpha-1}
\frac{(\lambda z^{-1}-1)^{\alpha-1}}{\lambda^{\alpha}} \frac{(2\lambda-1)^{m}}{(2\lambda+1)^{m+1}} d\lambda\cr
&
=&\frac{2}{\Gamma(\alpha)}
\int_{1}^{\infty} \frac{(u-1)^{\alpha-1}}{u^{\alpha}} \frac{(2uz-1)^{m}}{(2uz+1)^{m+1}} du\cr
&=&\frac{2(-1)^m}{2z\Gamma(\alpha)}
\int_{1}^{\infty} \frac{(u-1)^{\alpha-1}}{u^{\alpha}} \frac{(2(1/4z)-u)^{m}}{(2(1/4z)+u)^{m+1}}du\cr
&=&
{(-1)^m\over2z}\LL(\ell_{m,\a}) \left( \frac{1}{4z} \right).
\end{eqnarray*}
\end{proof}

\medskip
Let us now consider the space $\LL(\TT_2^{(\alpha)}(t^\a))$. Since $\LL$ is a
one-to-one mapping we define the norm $\Vert F\Vert_{\LL}:= \Vert f\Vert_{2,(\alpha)}$ on $\LL(\TT_2^{(\alpha)}(t^\a))$, where
$F=\LL(f)$, $f\in\TT_2^{(\alpha)}(t^\a)$. Clearly, $(\LL(\ell_{m,\a}))_{m,\a}$ is an orthonormal basis in
$\LL(\TT_2^{(\alpha)}(t^\a))$. Further, by Proposition \ref{Tau-properties} (iii)
and the classical Paley-Wiener theorem -that is, $H_2(\C^+)=\LL(L_2(\R^+))$- one has
$$
\LL(\TT_2^{(\beta)}(t^\beta))\hookrightarrow\LL(\TT_2^{(\alpha)}(t^\a))\hookrightarrow H_2(\C^+),\ \forall \beta\ge\alpha.
$$
Thus in particular we know that $\LL(\TT_2^{(\alpha)}(t^\a))$ is a RKHS.

\begin{proposition} \label{pro-nucleo} For $\a>0$, the reproducing kernel $Q_\a$ of
$\LL(\TT_2^{(\alpha)}(t^\a))$ is represented by the integral
$$
Q_{\a,w}(z):=Q_\a(z,w)= \int_{0}^{\infty} G_\a(z,r)\overline{G_a(w,r)}dr, \quad z, w\in\C^+,
$$
where
$$
G_\a(z,r)={\CC_\a(e_z)(r)\over\Gamma(\a+1)}
=
\frac{1}{r^{\a}}\int_{0}^{r}\frac{(r-u)^{\a-1}}{\G(\a)}e^{-zu}du, \quad z\in\C^+, r>0.
$$
\end{proposition}

\begin{proof}
Let $w\in\C^+$. By definition, there exists a unique $h_{\alpha,w}\in\TT_2^{(\a)}(t^\a)$ such that $Q_{\a,w}=\LL(h_{\alpha,w})$
in $\LL(\TT_2^{(\a)}(t^\a)$. Take $f\in\TT_2^{(\a)}(t^\a)$ and put $F=\LL f$ in $\LL(\TT_2^{(\a)}(t^\a))$.
Then by Hardy's inequality we have that the integral
$\displaystyle
\int_0^\infty s^\alpha\vert  W^\alpha f(s)\vert \left({1\over s^\alpha}
\int_0^s {(s-t)^{\alpha-1}e^{-(\Re w)t}dt}\right)ds$ is finite. Hence one can apply Fubini's theorem
(in the last-but-one equality of the following chain) to obtain
\begin{eqnarray*}
(\tau^\a W^\alpha f\vert \tau^\a W^\alpha h_{\alpha,w})_2&=&(f\vert h_{\alpha,w})_{2,(\alpha)}=(F\vert Q_{\a,w})_{\LL}=F(w)={\LL f}(w)\cr
&=&\int_0^\infty\int_t^\infty W^\alpha f(s){(s-t)^{\alpha-1}\over\Gamma(\alpha)} ds\ e^{-wt} dt\cr
&=&
\int_0^\infty s^\alpha W^\alpha f(s)\left({1\over s^\alpha}
\int_0^s {(s-t)^{\alpha-1}\over\Gamma(\alpha)} e^{-wt}dt\right)ds.
\end{eqnarray*}

Hence $Q_{\a,w}=\LL(h_{\alpha,w})$ with $h_{\alpha,w}$ such that
$$
s^\alpha W^\alpha h_{\alpha,w}(s)={1\over\Gamma(\a+1)}\CC_\a(e_{\overline{w}})(s),\ s>0;
$$
that is,
$$
h_{\alpha,w}(u)=\int_r^\infty {(u-r)^{\alpha-1}\over\Gamma(\a+1)}\CC_\a(e_{\overline{w}})(u){dr\over\Gamma(\alpha)}
={1\over\Gamma(\alpha+1)^2}\CC_\a^*(\CC_\a e_{\overline{w}})(u),\ u>0.
$$

Therefore
\begin{eqnarray*}
Q_\a(z,w)&=&{\LL(h_{\alpha,w})(z)\over\Gamma(\alpha+1)^2}
={(\CC_\a^*\CC_\a e_{\overline{w}}\vert e_{\overline{z}})_2\over\Gamma(\alpha+1)^2}
={(\CC_\a e_{\overline{w}}\vert \CC_\a e_{\overline{z}})_2\over\Gamma(\alpha+1)^2}
\cr
&=&\int_0^\infty
\left( \int_0^r {(r-s)^{\alpha-1}\over r^\alpha} {e^{-zs}ds\over \Gamma(\alpha)}\right)
\left(  \int_0^r {(r-t)^{\alpha-1}\over r^\alpha} {e^{-\overline{w}t} dt\over \Gamma(\alpha)}\right)\ dr,
\end{eqnarray*}
for every $z,w\in\C^+$, as we wanted to show.
\end{proof}

\begin{remark}\label{ISOgreen}
For $\alpha>1/2$, the above proposition is an immediate consequence of \cite[p. 82-83]{S1} when applied to
$\LL\colon\TT_2^{(\a)}(t^\a)\to H_2^{(\a)}(\C^+)$ since $\TT_2^{(\a)}(t^\a)$ is RKHS with kernel
$\int_{0}^{\infty} g_\a(s,r)\overline{g_a(t,r)}dr$, $s, t>0$, and $G_\a=\LL(g_\a)$. Recall that $\TT_2^{(\a)}(t^\a)$
is not a RKHS for $0<\alpha\le1/2$.
\end{remark}

\medskip
\begin{proof}({\it Theorem \ref{paleywiener})} Let $z,w\in\C^+$. By Proposition \ref{pro-nucleo},
\begin{eqnarray*}
Q_\a(z,w)&=&{1\over\Gamma(\a+1)^2}\int_{0}^{\infty}\CC_\a(e_z)(r)\CC_\a(e_{\overline{w}})(r)dr\cr
&=&\int_{0}^{\infty} \left( \int_{0}^r \frac{(r-s)^{\a-1}}{\G(\a)} \frac{e^{-zs}}{r^{\a}}ds \right)
\left( \int_{0}^r\frac{(r-t)^{\a-1}}{\G(\a)} \frac{e^{-\overline{w}t}}{r^{\a}} dt \right) dr\cr
&=& \int_{0}^{\infty} \int_{0}^{1} \frac{(1-x)^{\a-1}}{\G(\a)} \int_{0}^{1}
\frac{(1-y)^{\a-1}}{\G(\a)}e^{-zrx}e^{-\overline{w}ry}\ dx\ dy\ dr\cr
&=&\int_{0}^{1} \frac{(1-x)^{\a-1}}{\G(\a)} \int_{0}^{1}\frac{(1-y)^{\a-1}}{\G(\a)} \int_{0}^{\infty} e^{-r(xz+y\overline{w})}\ dr\ dy\ dx\cr
&=&\int_{0}^{1} \int_{0}^{1} \frac{(1-x)^{\a-1}}{\G(\a)} \frac{(1-y)^{\a-1}}{\G(\a)} \frac{1}{xz+y\overline{w}}\ dy\ dx.
\end{eqnarray*}

On the other hand,
\begin{eqnarray*}
K_\a(z,w)&=&\sum_{j=0}^\infty {\frak L}_{j,\alpha}(z)\overline{{\frak L}_{j,\alpha}(w)}={1\over 4z\overline w}
\sum_{j=0}^\infty{\mathcal L}(\ell_{j,\alpha})({1\over 4z})\overline{  {\mathcal L}(\ell_{j,\alpha})({1\over4w})}\\
&=&
{1\over 4z\overline w}\int_0^1\int_0^1{(1-x)^{\alpha-1}(1-y)^{\alpha-1}\over\Gamma(\alpha)^2}
{dx\ dy\over(x/4z)+(y/4\overline w)}\\
&=&\int_0^1\int_0^1{(1-x)^{\alpha-1}(1-y)^{\alpha-1}\over\Gamma(\alpha)^2}{dx\ dy\over x\overline {w}+yz}
=Q_\alpha(z,w),
\end{eqnarray*}
where the first and third equality hold true by \cite[Th. 1.8]{S1} and the second one follows from Lemma \ref{baseH2a}.

Since the reproducing kernels $K_\a$ and $Q_\a$ are equal the Hilbert spaces that they generate coincide, and the proof is over. \end{proof}

\medskip
\begin{remark}
\normalfont
Theorem \ref{paleywiener}(ii)  extends Theorem 3.3 of \cite{GMMS} to fractional $\a$. The proofs are rather different. In \cite{GMMS}, the identity $\LL(\TT_2^{(n)}(t^n))=H_2^{(n)}(\C^+)$, $n\in\N$, relies on the usage of Laguerre polynomials.
\end{remark}

%%%%%%%%%%%%%%%%%%%%%%%
%%%%%%%%%%%%%%%%%%%%%%%
\section{Estimating the kernel}\label{estkernel}
%%%%%%%%%%%%%%%%%%%%%
%%%%%%%%%%%%%%%%%%%%%

\medskip
We next proceed to estimate the norm of the kernel function $K_\a$. For $\a\ge1$, the calculation of such an estimate is the same as that one done for $\a=n\in\N$ in \cite{GMMS}. We include it here for the sake of completeness.

To begin with, one needs to know the value of the following standard integral.

\begin{lemma}\label{integ}
For $\theta\in ({-\pi/2}, {\pi/2})$,
$$
J(\theta):=\int_0^1{1\over t^2+1+2t\cos2\theta}\ dt={\vert \theta\vert\over\vert\sin 2\theta\vert}, \hbox{ if } \theta\not=0;
\quad J(0)=1/2.
$$
\end{lemma}

\begin{proof}
For $\theta=0$, one gets $J(0)=1/2$ very easily.
For $0<\vert\theta\vert<{\pi \over 2}$, we have
\begin{eqnarray*}
J(\theta) &= &\int_{0}^{1}\frac{dt}{(t+\cos(2\theta))^{2}+\sin^{2}(2\theta)}\cr
&=&
\frac{1}{\vert\sin(2\theta)\vert}
\int_{0}^{1/\vert\sin(2\theta)\vert}
\left[\left(r+\frac{\cos(2\theta)}{\vert\sin(2\theta)\vert}\right)^{2}+1\right]^{-1} dr\cr
&=&
\frac{1}{\vert\sin(2\theta)\vert}\int_{\cos(2\theta)/\vert\sin(2\theta)\vert}^{\cos\theta/\vert\sin\theta\vert}\frac{du}{u^{2}+1}
= {\vert\theta\vert\over\vert\sin(2\theta)\vert},
\end{eqnarray*}
as we wanted to show.
\end{proof}

\begin{theorem}\label{estimation}
Let $\a>0$. Then, for every $z=\vert z\vert e^{i\theta}\in\C^+$,
\item{\rm{(i)}} for $\a\ge1$,
$$
{1\over (2\alpha-1)\Gamma(\alpha)^2}\ {1\over\vert z\vert}\le\Vert K_\alpha(\cdot,z)\Vert_{2,(\a)}^2
\le{\pi\over \alpha\Gamma(\alpha)^2}\ {1\over\vert z\vert};
$$
\item{\rm{(ii)}} for $1/2<\a<1$,
$$
{1\over \Gamma(\alpha)^2}\ {1\over\vert z\vert}\le\Vert K_\alpha(\cdot,z)\Vert_{2,(\a)}^2
\le{\pi\over (2\alpha-1)\Gamma(\alpha)^2}\ {1\over\vert z\vert};
$$
\item{\rm{(iii)}} for $0<\a\le1/2$,
$$
{1\over\Gamma(\a)^2}{1\over\vert z\vert}\le\Vert K_\alpha(\cdot,z)\Vert_{2,(\a)}^2\le{2\over\Gamma(\a+1)^2}{1\over\vert z\vert},
\hbox{ if } \vert\theta\vert\le\pi/4,
$$
and
$$
{1\over\Gamma(\a)^2}{1\over\vert z\vert}\le\Vert K_\alpha(\cdot,z)\Vert_{2,(\a)}^2\le{2\over\Gamma(\a+1)^2}{1\over\Re z},
\hbox{ if } \pi/4<\vert\theta\vert<\pi/2.
$$
\end{theorem}

\begin{proof}
Let $z=\vert z\vert e^{i \theta}\in\CCma$ with $\theta \in (-\pi/2, \pi/2)$. For every $\a>0$ one has
\begin{eqnarray*}
\Vert K_{\a,z}\Vert_{2,(\a)}^2&=&K_\a(z,z)
=\int_{0}^{1}\int_{0}^{1}\frac{(1-y)^{\a-1}}
{\Gamma(\a)}\frac{(1-x)^{\a-1}}{\Gamma(\a)}\frac{1}{x{\overline z}+yz}dxdy\cr
&=&
\frac{1}{|z|}\int_{0}^{1}\int_{0}^{1}\frac{(1-y)^{\a-1}}{\Gamma(\a)}\frac{(1-x)^{\a-1}}{\Gamma(\a)}\frac{(x+y)\cos\theta}{x^{2}+y^{2}+2xy\cos(2\theta)}dxdy,
\end{eqnarray*}
since $\Im(K_{n,z}(z))=0$ (use symmetry in the imaginary part of the integral).
Thus using symmetry again (with respect to the diagonal in $(0,1)\times(0,1)$) and then the change of variables
$x=yt$ one obtains
\begin{eqnarray*}
K_\a(z,z) &=&
\frac{2\cos\theta}{\Gamma(\a)^2|z|}
\int_{0}^{1}\int_{0}^{y}\frac{(1-x)^{\a-1}(1-y)^{\a-1}(x+y)}{x^{2}+y^{2}+2xy\cos(2\theta)}dxdy\cr
&\buildrel{(*)}\over=&\frac{2\cos\theta}{\Gamma(\a)^2|z|}
\int_{0}^{1}\left(\int_{0}^{1}(1-yt)^{\a-1}(1-y)^{\a-1}dy\right){(1+t)\over t^{2}+1+2t\cos(2\theta)}dt.
\end{eqnarray*}

(i) If $\a\ge1$ then
$$
\int_{0}^{1}(1-yt)^{\a-1}(1-y)^{\a-1}dy\le\int_{0}^{1}(1-y)^{\a-1}dy={1\over\a}
$$
and therefore, by equality $(*)$,
$$
K_\a(z,z)\le
\frac{4\cos\theta}{\a\Gamma(\a)^2|z|}J(\theta)
={2\over\a\Gamma(\a)^2}{\vert\theta\vert\over\vert\sin\theta\vert}{1\over\vert z\vert}
\le{\pi\over\a\Gamma(\a)^2}{1\over\vert z\vert}.
$$

For the lower estimate, one has $(1-yt)^{\a-1}\ge (1-y)^{\a-1}$ for $0<t,y<1$ and so
\begin{eqnarray*}
K_\a(z,z)
&\ge&\frac{2\cos\theta}{\Gamma(\a)^2|z|}\int_{0}^{1}\int_{0}^{1}\frac{(1-y)^{2\a-2}}{t^{2}+1+2t\cos(2\theta)}dtdy\cr
&=&
\frac{2}{\Gamma(\a)^2(2\a-1)}{\vert\theta\vert\over\vert\sin\theta\vert}{1\over\vert z\vert}\ge
\frac{1}{\Gamma(\a)^2(2\a-1)}{1\over\vert z\vert}.
\end{eqnarray*}

(ii) If $1/2<\a<1$, using the Cauchy-Schwarz inequality, one obtains
\begin{eqnarray*}
\int_{0}^{1}(1-yt)^{\a-1}(1-y)^{\a-1}dy&\le&{1\over2\a-1}\left({1-(1-t)^{2\a-1}\over t}\right)^{1/2}\cr
&\le&
\left(\sup_{0<t<1}{1-(1-t)^{2\a-1}\over(2\a-1)t}\right)={1\over2\a-1}.
\end{eqnarray*}
and then,
by $(*)$ prior to the proof of  (i),
\begin{eqnarray*}
K_\a(z,z)&\le&
\frac{4\cos\theta}{\Gamma(\a)^2(2\a-1)}J(\theta){1\over\vert z\vert}
\le\frac{\pi}{\Gamma(\a)^2(2\a-1)}{1\over\vert z\vert}.
\end{eqnarray*}

As for the lower estimate, from $(*)$
one gets
\begin{eqnarray*}
K_\a(z,z)&\ge&\frac{2\cos\theta}{\Gamma(\a)^2|z|}\int_{0}^{1}\int_{0}^{1}\frac{dy\ dt}{t^{2}+1+2t\cos(2\theta)}\cr
&=&\frac{2\cos\theta}{\Gamma(\a)^2|z|}J(\theta)\ge\frac{1}{\Gamma(\a)^2}{1\over|z|}.
\end{eqnarray*}

(iii) Assume $0<\a\le1/2$. We now deal first with the lower estimate. Then, as in (ii),
$$
K_\a(z,z)\ge\frac{2\cos\theta}{\Gamma(\a)^2|z|}J(\theta)=\frac{1}{\Gamma(\a)^2}{1\over|z|}.
$$

For the upper estimate,
by $(*)$,
\begin{eqnarray*}
K_\a(z,z)&\le&\frac{4\cos\theta}{\Gamma(\a)^2|z|}\left(\int_{0}^{1} (1-y)^{\a-1}dy\right)
\left(\int_0^1{(1-t)^{\a-1}\over t^{2}+1+2t\cos(2\theta)}dt\right)\cr
&=&\frac{4\cos\theta}{\a\Gamma(\a)^2|z|}\int_0^1{(1-t)^{\a-1}\over t^{2}+1+2t\cos(2\theta)}dt.
\end{eqnarray*}

Now is the moment to notice that $\cos(2\theta)\ge0$ if $\vert\theta\vert\le\pi/4$ whereas $\cos(2\theta)<0$ if $\vert\theta\vert>\pi/4$.

For $\vert\theta\vert\le\pi/4$,
$$
K_\a(z,z)\le\frac{4\cos\theta}{\a\Gamma(\a)^2|z|}\int_0^1{(1-t)^{\a-1}\over t^{2}+1}dt
\le\frac{2\cos\theta}{\Gamma(\a+1)^2|z|}
\le\frac{2}{\Gamma(\a+1)^2|z|}.
$$

For $\vert\theta\vert>\pi/4$,
\begin{eqnarray*}
K_\a(z,z)&\le&\frac{4\cos\theta}{\a\Gamma(\a)^2|z|}\int_0^1{(1-t)^{\a-1}\over (t+\cos(2\theta))^2+\sin^2(2\theta)}dt\cr
&\le&\frac{4\cos\theta}{\a\Gamma(\a)^2|z|}\int_0^1{(1-t)^{\a-1}\over \sin^2(2\theta)}dt\cr
&=&{1(1/\sin^2\theta)\over\Gamma(\a+1)^2\cos\theta}{1\over\vert z\vert}
\le{2\over\Gamma(a+1)^2}{1\over \Re z}.
\end{eqnarray*}

The proof is over.
\end{proof}

In the above theorem, parts (i) and (ii), we have shown that
$$
\Vert K_\alpha(\cdot,z)\Vert_{2,(\alpha)}\simeq {1\over\sqrt{\vert z\vert}}, \ z\in\C^+,
$$
up to constants only depending on $\a$, for $\a>1/2$.

\medskip
QUESTION.- In Theorem \ref{estimation} (iii), so for $0<\a\le1/2$,
\begin{itemize}
\item[(i)] Is $\Vert K_\alpha(\cdot,z)\Vert_{2,(\alpha)}\simeq |z|^{-1/2}, \ z\in\C^+$ ?
\item[(ii)] Is $\Vert K_\alpha(\cdot,z)\Vert_{2,(\alpha)}\simeq (\Re z)^{-1/2}, \ z\in\C^+ $ ?
\item[(iii)] Or maybe, does not any of the above hold true ?
\end{itemize}

%%%%%%%%%%%%%%%%%%%%%%%
%%%%%%%%%%%%%%%%%%%%%%%
\section{Averaging Brownian motion}\label{brownaverage}
%%%%%%%%%%%%%%%%%%%%%
%%%%%%%%%%%%%%%%%%%%%

One can also introduce a Hilbert space of absolutely continuous functions of fractional order, using the Ces\`aro-Hardy operator
$\CC_\a$ in a similar manner to which we have done using $\CC_\a^*$ above. We focus on the case $p=2$.

Set $\TT_{(\a)}^2:=\CC_\a(L_2(\R^+))$. Since $\CC_\a$ is injective by Tichsmarsh's theorem we are allowed to define the norm
$\Vert f\Vert_{(\a),2}:=\Gamma(\a+1)\Vert (\CC_\a)^{-1}f\Vert_2$, for all $f\in\TT_{(\a)}^2$, so that
$\TT_{(\a)}^2$ becomes a Hilbert space with inner product
$$
(f\vert g)_{(\a),2}=\Gamma(\a+1)^2\int_0^\infty((\CC_\a)^{-1}f)(s)\overline{((\CC_\a)^{-1}g)(s)}ds,\ f,g\in\TT_{(\a)}^2,
$$
such that $\TT_{(\a)}^2\hookrightarrow L_2(\R^+)$. Moreover, the formula
$$
f(t)={\a\over t^\a}\int_0^t(t-s)^{\a-1}((\CC_\a)^{-1}f)(s) ds
$$
entails that $\TT_{(\a)}^2$ is a RKHS for $\a>1/2$.

Note that, in the notation of Section \ref{brownian}, $\DD^\a(\tau_\a f)=\Gamma(\a+1)\CC_\a^{-1}(f)$.

Set $H_{(\a)}^2(\C^+):=\LL(\TT_{(\a)}^2)$ endowed with the inner product
$$
(\LL f\vert\LL g)_{(\a),2}:=( f\vert g)_{(\a),2}=(\DD^\a(\tau_\a f)\vert\DD^\a(\tau_\a g))_2.
$$

Since $H_{(\a)}^2(\C^+)\hookrightarrow H_2(\C^+)$ the space $H_{(\a)}^2(\C^+)$ is a RKHS for every $\a>0$. Let $R_\a$ its
reproducing kernel. There is a unique $\phi_{w,\a}\in\TT_{(\a)}^2$ such that $\LL(\phi_{w,\a})=R_{\a,w}$. Then, for $w\in\C^+$
and every
$f\in\TT_{(\a)}^2$,
\begin{eqnarray*}
(\DD^\a(\tau_\a f)\vert\DD^\a(\tau_\a \phi_{w,\alpha})_2&=&(\LL f\vert R_{\a,w})_{(\a),2}=\LL f(w)\cr
&=&(\CC_\a\CC_\a^{-1}f)\vert e_{\overline{w}})_2
=(\CC_\a^{-1}f\vert\CC_\a^*( e_{\overline{w}}))_2,
\end{eqnarray*}
and therefore we derive
$\CC_\a^{-1}\phi_{w,\a}=\Gamma(\a+1)^{-2}\ \CC_\a^*( e_{\overline{w}})$, whence
$$
\phi_{w,\a}={1\over\Gamma(\a+1)^2}\CC_\a\CC_\a^*( e_{\overline{w}})=
{1\over\Gamma(\a+1)^2}\CC_\a^*\CC_\a( e_{\overline{w}}),
$$
by Corollary \ref{conmuta}.

In other words, we have found that $\phi_{w,\a}=h_{\a,w}$ where $h_{\a,w}$ is the function obtained prior to
Remark \ref{ISOgreen}. This implies that
$R_\a(z,w)=Q_\a(z,w)=K_\a(z,w)$ for $z, w\in\C^+$,
which is to say that the  spaces $H_{(\a)}^2(\C^+)$ and $H^{(\a)}_2(\C^+)$ are the same
(and so are $\TT_2^{(\a)}(t^\a)$ and $\TT^2_{(\a)}$).
In particular, we have got the equivalence of the norms
$\Vert\DD^\a(\tau_\a f)\Vert_2$ and
$\Vert \tau_\a W^\a f\Vert_2$ for $f\in\TT^2_{(\a)}$,
which can be considered as an extension of  \cite[Prop. 2.6]{GMMS} to fractional derivatives.

\begin{remark}
\normalfont
Note that the argument followed prior to this remark allows us to identify the spaces $\LL(\CC_\a(L_2(\R^+)))$ and
$\LL(\TT_2^{(\a)}(t^\a))$
--previous identification of $\phi_{w,\a}$ and $h_{\a,w}$-- {\it independently} of
Theorem \ref{paleywiener}. In fact, part (ii) of such a theorem is easily obtained from the equality
$\LL(\CC_\a(L_2(\R^+)))=\LL(\TT_\a^{(\a)}(t^\a))$ using Corollary \ref{laplacesaro}:
\begin{eqnarray*}
\LL(\CC_\a(L_2(\R^+)))&=&(\LL\circ\CC_\a)(L_2(\R^+))=({\frak C}_\a\circ\LL)(L_2(\R^+))\cr
&=&{\frak C}_\a(\LL(L_2(\R^+)))={\frak C}_\a(H_2(\C^+))=H^{(\a)}_2(\C^+)
\end{eqnarray*}
where we have used that corollary in the second equality and the classical Paley-Wiener theorem in the last but one equality.

We have however chosen to present our first approach to Theorem 6.2 to give a more complete description of the space
$H_2^{(\a)}(\C^+)$, by showing in particular a nice basis in it. Recall in passing that having on hand suitable bases in RKHS of the fBm's may be helpful to get representations of Gaussian processes, see \cite[p. 3]{Hu}.

On the other hand, self-similarity is an important topic related with fBm's where in particular Hardy spaces on
$\C^+$ of fractional derivatives of Laplace transforms are of interest, see \cite[p. 277]{M}.
In this respect, notice that the kernels (covariances) $k_\a$ and $K_\a$ of preceding sections satisfy, for $\lambda>0$, $s,t>0$ and
$z,w\in\C^+$, the rules
$$
k_\a(\lambda s,\lambda t)=\lambda^{-1}k_\a(s,t); \quad  K_\a(\lambda z,\lambda w)=\lambda^{-1}K_\a(z,w)
$$
{\it independently of} $\a>0$.
\end{remark}

{\begin{remark}
\normalfont The identification between the two spaces $\CC_\alpha(L_2(\R^+))$ and $\CC_\alpha^*(L_2(\R^+))$ very much depends on the existence of reproducing  kernels as well as of the isomorphisms provided by the Laplace transform in the Hilbertian setting. For general 
$1\le p(\not=2)\le \infty$ these tools fail and it is not clear how to substitute them with other suitable to characterize those range spaces. We leave this question open.
\end{remark}

In conclusion, taking into account the properties of spaces $\TT_2^{(\a)}(t^\a)$ pointed out in previous sections and the simple but interesting description of the space $H^{(\a)}(\C^+)$ of Laplace transforms of the elements in $\TT_2^{(\a)}(t^\a)$, we wonder if operating with averages of fractal Brownian processes
$t^{-\a}\int_0^t(t-s)^{\a-1} B_s ds$ or $t^{-\a}\int_0^t(t-s)^{\a-1} h(s) ds$, $h\in L_2(\R^+)$, could be helpful in this setting, at least from an operational viewpoint.

\medskip

ACKNOWLEDGEMENT.- We wish to thank the referee for comments and remarks which have contributed to improve the presentation of this paper. 

%%%%%%%%%%%%%%%
%%%%%%%%%%%%%%%

\end{document}